\documentclass[10pt,a4paper]{article}
\usepackage{amssymb,latexsym,times,standard,graphicx,amsmath}
\usepackage{pdfsync}
\usepackage{mathtools,bm,mathrsfs}
\usepackage[title]{appendix}
\textheight 
            675pt
\textwidth 
           460pt

\newcommand{\st}{\ \big| \ }

\newcommand{\dk}{\, \mathrm{d}k}
\newcommand{\dtildek}{\, \mathrm{d}\tilde{k}}

\newcommand{\dtilder}{\, \mathrm{d}\tilde{r}}
\newcommand{\dr}{\, \mathrm{d}r}

\newcommand{\dz}{\, \mathrm{d}z}

\newcommand{\bfB}{{\bm{B}}}

\newcommand{\bfH}{{\bm{H}}}

\newcommand{\bfM}{{\bm{M}}}

\newcommand{\bfv}{{\bm{v}}}

\newcommand{\FF}{\mathcal F}
\newcommand{\GG}{\mathcal G}
\newcommand{\HH}{\mathcal H}

\newcommand{\KK}{\mathcal K}
\newcommand{\LL}{\mathcal L}

\newcommand{\NN}{\mathcal N}
\newcommand{\OO}{\mathcal O}

\newcommand{\RR}{\mathcal R}

\newcommand{\WW}{\mathcal W}
\newcommand{\XX}{\mathcal X}

\newcommand{\ZZ}{\mathcal Z}

\newcommand{\e}{\mathrm{e}}
\renewcommand{\i}{\mathrm{i}}

\DeclareMathOperator{\sech}{sech}
\DeclareMathOperator{\re}{Re}
\DeclareMathOperator{\im}{Im}

\newcommand{\qed}{\hfill$\Box$\bigskip}

\newcommand{\alignqed}{\tag*{$\Box$}}

\newcommand{\nn}{|{\mskip-2mu}|{\mskip-2mu}|}

\newtheorem{theorem}{Theorem}[section]
\newtheorem{lemma}[theorem]{Lemma}

\newtheorem{proposition}[theorem]{Proposition}
\newtheorem{corollary}[theorem]{Corollary}
\newtheorem{remark}[theorem]{Remark}
\newtheorem{remarks}[theorem]{Remarks}

\begin{document}
\newcounter{count}

\title{An existence theory for solitary waves on a ferrofluid jet}

\author{M. D. Groves\thanks{FR Mathematik, Universit\"{a}t des Saarlandes, Postfach 151150, 66041 Saarbr\"{u}cken, Germany}
\and
D. Nilsson\thanks{Department of Mathematics, Linnaeus University, V\"{a}xj\"{o}, Sweden}
\and
L. Sch\"{u}tz\footnotemark[1]}
\date{}
\maketitle

\begin{abstract}
We discuss axisymmetric solitary waves on the surface of an otherwise cylindrical ferrofluid jet surrounding a stationary metal rod.
The ferrofluid, which is governed by a general (nonlinear) magnetisation law, is subject to an azimuthal magnetic field generated by an electric current flowing along the
rod. We treat the governing equations using a modification of the Zakharov-Craig-Sulem formulation for water waves, reducing the problem to a
single nonlocal equation for the free-surface elevation variable $\eta$. The nonlocality in the equation takes the form of a Dirichlet-Neumann operator whose analyticity
(in standard function spaces) is demonstrated by studying its defining boundary-value problem in newly introduced Sobolev spaces for radial functions.\
Using rudimentary fixed-point arguments and
Fourier analysis we rigorously reduce the equation for $\eta$ to a perturbation of a Korteweg-de Vries equation (for strong surface tension) or a nonlinear Schr\"{o}dinger equation
(for weak surface tension), both of which have nondegenerate explicit solitary-wave solutions. The existence theory is completed using an appropriate version of the implicit-function theorem.
\end{abstract}

\section{Introduction}

\subsection{The hydrodynamic problem}
\begin{figure}[h]
\centering
\includegraphics{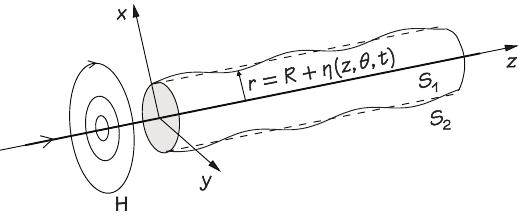}
\caption{Waves on the surface of a ferrofluid jet surrounding a current-carrying wire} \label{Configuration}
\end{figure}

We consider the inviscid, incompressible and irrotational flow of a ferrofluid of unit density in the region
$$S_1=\{0<r<R+\eta(\theta,z,t)\}$$
bounded by the free surface $\{r=R+\eta(\theta,z,t)\}$ and a current-carrying wire at $\{r=0\}$ (see Figure \ref{Configuration}). Here $(r,\theta,z)$
are the usual cylindrical polar coordinates, $t$ is time, $R$ is a positive constant which represents the radius of the jet without any current flow,
and $\eta$ is a function of $(\theta,z,t)$. The magnetic field generated by the wire is static and the region
$$S_2 = \{r>R+\eta(\theta,z,t)\}$$
is assumed to be a vacuum. The irrotational magnetic and solenoidal induction fields in $S_1$ and $S_2$ are denoted by respectively $\bfH_1$, $\bfB_1$ and $\bfH_2$, $\bfB_2$,
while the irrotational, solenoidal velocity field of the fluid in $S_1$ is denoted by $\bfv$. The interdependence between the fields is given by the formulae
$$\bfB_1 = \mu_0(\bfH_1 + \bfM_1(\bfH_1)), \qquad \bfB_2=\mu_0 \bfH_2,$$
where $\mu_0$ is the magnetic permeability of free space,
$$\bfM_1(\bfH_1) = m_1(|\bfH_1|)\frac{\bfH_1}{|\bfH_1|}$$
is the given magnetic intensity of the ferrofluid and $m_1({\bfH_1})$ is a nonnegative function.

The ferrohydrodynamic problem was formulated in terms of magnetic potential functions $\psi_1$, $\psi_2$ and
a velocity potential $\phi$ such that
$$\bfH_1 = - \nabla \psi_1, \quad \bfH_2 = -\nabla \psi_2, \quad \bfv=\nabla \phi$$
by Groves \& Nilsson \cite[\S2]{GrovesNilsson18} following the theory given by Rosensweig \cite[\S\S5.1--5.2]{Rosensweig}. The governing equations are
\begin{align}
\nabla \cdot (\mu(|\nabla \psi_1|\nabla \psi_1)) =0, \quad\qquad & 0<r<R+\eta(\theta,z,t), \label{Body 1}\\
\Delta \psi_2 =0, \quad\qquad & r>R+\eta(\theta,z,t), \label{Body 2}\\
\Delta \phi =0, \quad\qquad & 0<r<R+\eta(\theta,z,t), \label{Body 3}
\end{align}
where
$$\mu(s)=1+\frac{m(s)}{s},$$
with boundary conditions
\begin{align}
\psi_2-\psi_1 &=0, \label{BC 1}\\
\psi_{2n}-\mu(|\nabla\psi_1|) \psi_{1n} &=0, \label{BC 2}\\
-\eta_t + \phi_r - \frac{1}{r^2}\phi_\theta\eta_\theta - \phi_z \eta_z &=0, \label{BC 3} \\
\phi_t+\tfrac{1}{2}|\nabla\phi|^2-\mu_0\nu(|\nabla\psi_1|)+2\sigma \kappa - \tfrac{1}{2}\mu_0 (\mu(|\nabla\psi_1|)-1)^2 &=c_0 \label{BC 4}
\end{align}
at $r=R+\eta(\theta,z,t)$, where $2\kappa$
is the mean curvature of the surface, $\sigma$ is the coefficient of surface tension and
$c_0$ is a constant arising from integration of the (magnetic) Euler equation. Equations \eqref{Body 1}--\eqref{Body 3} state
that $\bfB_1$, $\bfB_2$ and $\bfv$ are solenoidal, equations \eqref{BC 1}, \eqref{BC 2} state that the magnetic and induction fields
are continuous at the surface, while equation \eqref{BC 3} is the the hydrokinematic boundary condition that fluid particles on the surface remain there and
equation \eqref{BC 4} is the hydrodynamic boundary condition which balances the forces at the surface.

The constant $c_0$ is selected so that
$${\bf H}_1=\dfrac{J}{2\pi r}{\bf e}_\theta, \qquad{\bf H}_2=\dfrac{J}{2\pi r}{\bf e}_\theta, \qquad {\bf v}={\bf 0}, \qquad \eta=0$$
(that is $\psi_1=\psi_2=-J\theta/2\pi$, $\phi=0$, $\eta=0$)
is a solution to the above equations
(corresponding to a
uniform magnetic field and a circular cylindrical jet with radius $R$); we therefore set
$c_0 = -\mu_0\nu(J/2\pi r)+\sigma/R$.
Seeking axisymmetric waves for which $\eta$ and $\phi$ are independent of $\theta$, one finds that
$\psi_1=\psi_2=-J\theta/2\pi$, so that the hydrodynamic problem decouples from the magnetic problem
and is given by
$$
\phi_{rr}+\frac{1}{r}\phi_r + \phi_{zz} =0, \quad\qquad 0<r<R+\eta(z,t)
$$
and
\begin{align*}
& -\eta_t+\phi_r-\phi_z\eta_z = 0, \\
& \phi_t+\frac{1}{2}(\phi_r^2+\phi_z^2) - \mu_0\nu\left(\frac{J}{2\pi(R+\eta)}\right)
\\
& \hspace{1cm}\mbox{}+\mu_0\nu\left(\frac{J}{2\pi R}\right)+\frac{\sigma}{(R+\eta)(1+\eta_z^2)^{1/2}} - \frac{\sigma \eta_{zz}}{(1+\eta_z^2)^{3/2}}-\frac{\sigma}{R}=0
\end{align*}
at $r=R+\eta(z,t)$, where we have used the formula
$$2\kappa = \frac{- (R+\eta)^2(1+\eta_z^2)+(R+\eta)^3\eta_{zz}}
{(R+\eta)^{3/2}(1+\eta_z^2)^{3/2}}.$$

This initial-value problem has been studied by Wang \& Yang \cite{WangYang19}, but here we concentrate upon travelling waves.
Introducing dimensionless variables
\[
 (\hat z,\hat r):=\frac{1}{R}(z,r),\qquad \hat{t}=\frac{\sigma^{1/2}}{R^{3/2}}t,
 \qquad \hat \phi:= \frac{1}{(\sigma R)^{1/2}}\phi,\qquad
\hat\eta:= \frac{1}{R}\eta
\]
and functions
\[
\hat{m}_1(s):=\frac{2\pi R}{J \chi}m_1\left(\frac{J}{2\pi R}s\right), \qquad
\hat\nu(s):= \frac{4\pi^2 R^2}{J^2\chi}\nu\left(\frac{J}{2\pi R}s\right),
\]
where $\chi=(2\pi R/J)m_1(J/2\pi R)$ and
$\hat{m}(1)=\hat{\nu}^\prime(1)=1$, and looking for travelling-wave solutions of the form
$$\phi(r,z,t)=\phi(r,z-ct), \qquad \eta(z,t)=\eta(z-ct),$$
we arrive at the equations
\begin{equation}
\phi_{rr}+\frac{1}{r}\phi_r + \phi_{zz} =0, \qquad 0<r<1+\eta(z,t), \label{Unflattened 1}
\end{equation}
and
\begin{align}
& c\eta_z+\phi_r-\phi_z\eta_z = 0, \label{Unflattened 2}\\
& -c\phi_z+\frac{1}{2}(\phi_r^2+\phi_z^2)-\gamma\left(\nu\left(\frac{1}{1+\eta}\right)-\nu(1)\right)
+\left(\frac{1}{(1+\eta)(1+\eta_z^2)^{1/2}} - \frac{\eta_{zz}}{(1+\eta_z^2)^{3/2}}-1\right)=0 \label{Unflattened 3}
\end{align}
at $r=1+\eta(z,t)$, where 
$$\gamma=\frac{\mu_0 J^2\chi}{4\pi^2\sigma R^2}.$$
Solitary waves are nontrivial solutions to \eqref{Unflattened 1}--\eqref{Unflattened 3} which are evanescent as $|z| \rightarrow \infty$.

\subsection{The main results}

We treat equations \eqref{Unflattened 1}--\eqref{Unflattened 3} using a modification of the Zakharov-Craig-Sulem formulation for water waves
(Zakharov \cite{Zakharov68}, Craig \& Sulem \cite{CraigSulem93}), thus reducing the problem to a single non-local equation for $\eta$ by introducing a
\emph{Dirichlet-Neumann operator} informally defined as follows (see Xu \& Wang \cite{XuWang25} for a similar approach for the time-dependent problem
and Blyth \& Parau \cite{BlythParau19} for an alternative non-local reformulation). Fix $\Phi=\Phi(z)$, let $\phi$ be the unique solution of the
Dirichlet boundary-value problem
\begin{align}
\phi_{rr}+\frac{1}{r}\phi_r + \phi_{zz} &=0, \quad\qquad 0<r<1+\eta, \label{DNO defn 1}\\
\phi &= \Phi \quad\qquad r=1+\eta, \label{DNO defn 2}
\end{align}
and define
\begin{align*}
G(\eta)\Phi & := (1+\eta)(1+\eta_z^2)^{1/2}\frac{\partial \phi}{\partial n}\Big|_{r=1+\eta} \\
& = (1+\eta)(\phi_r - \eta_z\phi_z)\big|_{r=1+\eta}.
\end{align*}
Equations \eqref{Unflattened 2} and \eqref{Unflattened 3} can be rewritten as
\begin{equation}
c \eta_z + \frac{G(\eta)\Phi}{1+\eta}=0 \label{Reform BC 1}
\end{equation}
and
\begin{align}
-c\Phi_z &+ \tfrac{1}{2}\Phi_z^2 - \frac{1}{2(1+\eta_z^2)}\left(\eta_z\Phi_z + \frac{G(\eta)\Phi}{1+\eta}\right)^{\!\!2} \nonumber \\
& \!\!\!-\gamma\left(\nu\left(\frac{1}{1+\eta}\right)-\nu(1)\right)
+\left(\frac{1}{(1+\eta)(1+\eta_z^2)^{1/2}} - \frac{\eta_{zz}}{(1+\eta_z^2)^{3/2}}-1\right)=0, \label{Reform BC 2}
\end{align}
and by substituting $\Phi=-cG(\eta)^{-1}(\eta_z+\eta\eta_z)$ from \eqref{Reform BC 1} into \eqref{Reform BC 2}, we arrive at
\begin{equation}
\KK(\eta)-c^2\LL(\eta) =0, \label{GZCS}
\end{equation}
where
\begin{align}
\KK(\eta) &= -\gamma\left(\nu\left(\frac{1}{1+\eta}\right)-\nu(1)\right)
+\left(\frac{1}{(1+\eta)(1+\eta_z^2)^{1/2}} - \frac{\eta_{zz}}{(1+\eta_z^2)^{3/2}}-1\right), \label{Defn of KK}\\
\LL(\eta) & = -\tfrac{1}{2}(K(\eta)\eta+\tfrac{1}{2}K(\eta)\eta^2)^2 + \frac{1}{2(1+\eta_z^2)}(\eta_z -\eta_z K(\eta)\eta-\tfrac{1}{2}\eta_zK(\eta)\eta^2)^2 \nonumber \\
& \qquad \mbox{}+K(\eta)\eta+\tfrac{1}{2}K(\eta)\eta^2 \label{Defn of LL}
\end{align}
and
\begin{equation}
K(\eta)\xi = -(G(\eta)^{-1}\xi_z)_z. \label{Defn of K}
\end{equation}
Equation \eqref{GZCS} is equivalent to \eqref{Unflattened 1}--\eqref{Unflattened 3}; the velocity potential is recovered by setting
$\Phi=-cG(\eta)^{-1}(\eta_z+\eta\eta_z)$ and solving \eqref{DNO defn 1}, \eqref{DNO defn 2}.

Our task is therefore to find nontrivial solutions to \eqref{GZCS} which satisfy $\eta(z) \rightarrow 0$ as $z \rightarrow \pm \infty$, and we prove the following results.

\begin{theorem} \label{KdV thm}
Suppose that $1<\gamma<9$ and $c^2=c_0^2(1-\varepsilon^2)$. For each sufficiently small value of $\varepsilon>0$ there exists a symmetric
\underline{Korteweg-de Vries solitary-wave solution} of \eqref{GZCS} which satisfies
$$\eta(z)= \varepsilon^2 \zeta_\mathrm{KdV}(\varepsilon z) + o(\varepsilon^2)$$
uniformly over $z \in {\mathbb R}$, where
\begin{equation}
\zeta_\mathrm{KdV}(Z)=-\frac{3}{2d_0}\sech^2 \left(2 \left(\frac{c_0^2}{9-\gamma }\right)^{1/2}Z\right) \label{Explicit KdV}
\end{equation}
and
$$d_0 = \frac{1}{2c_0^2}\left(\frac{3}{2}\gamma - \frac{1}{2}\gamma\nu^{\prime\prime}(1)-\frac{3}{2}\right), \qquad c_0^2 = \frac{1}{2}(\gamma-1).$$
\end{theorem}

\begin{figure}[h]
\centering
\includegraphics[scale=0.3]{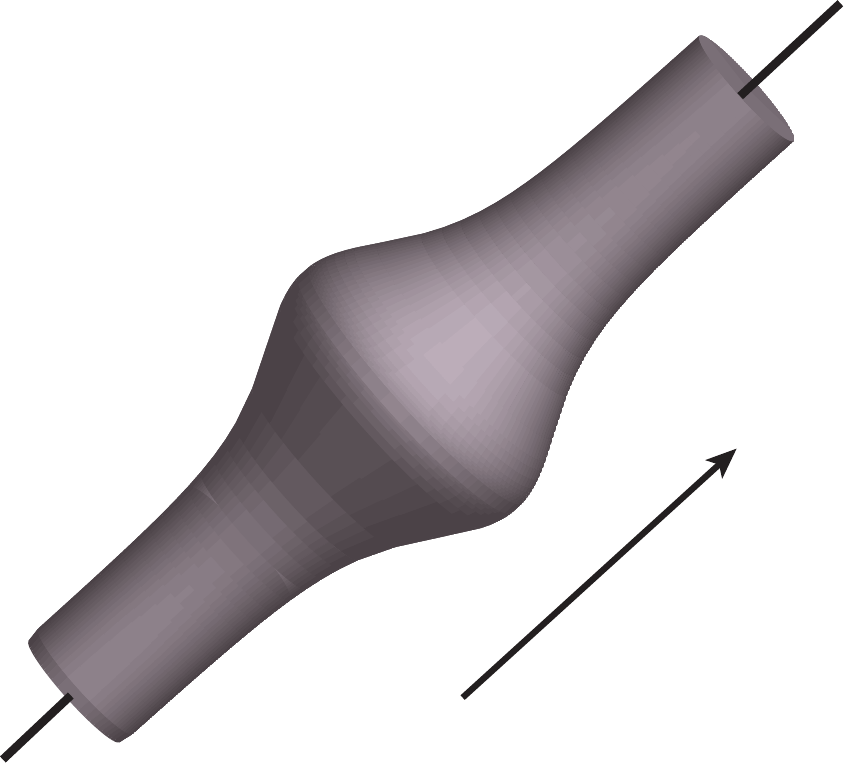}\hspace{1in}\includegraphics[scale=0.3]{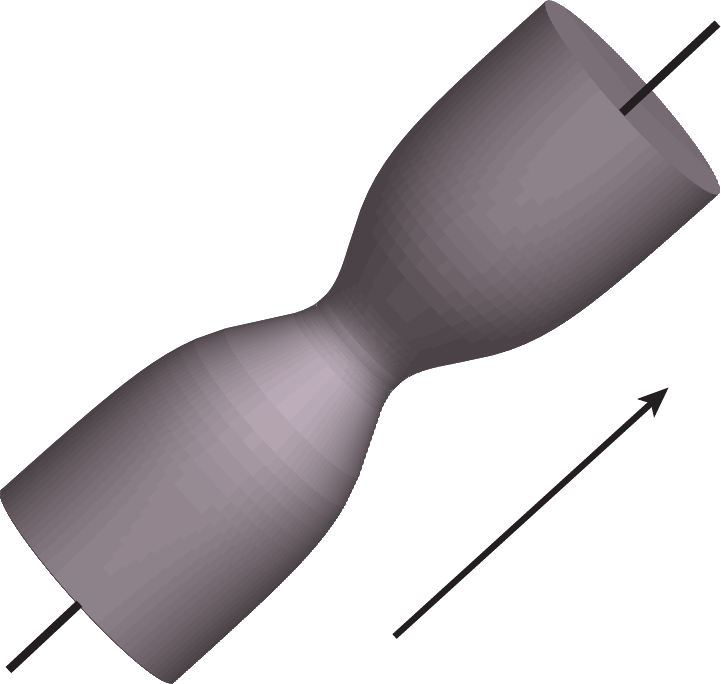}
\caption{Korteweg-de Vries solitary waves of elevation (left) and of depression (right) depending on the sign of $d_0$}
\end{figure}

\begin{theorem} \label{NLS thm}
Suppose that $\omega>0$,
$$\gamma=1-\omega^2 + \frac{2\omega f(\omega)}{f^\prime(\omega)}, \qquad c_0^2 = \frac{2\omega}{f^\prime(\omega)}, \qquad f(\omega)=\frac{\omega I_0(\omega)}{I_1(\omega)},$$
where $I_\nu$ is the modified Bessel function of the first kind and order $\nu$, and $c^2=c_0^2(1-\varepsilon^2)$.
For each sufficiently small value of $\varepsilon>0$ there exist two symmetric
\underline{nonlinear Schr\"{o}dinger solitary-wave solutions} of \eqref{GZCS} which satisfy
\begin{equation}
\eta(z) = \pm\varepsilon \zeta_\mathrm{NLS}(\varepsilon z)\cos(\omega z) + o(\varepsilon) \label{eta from zeta}
\end{equation}
uniformly over $z \in {\mathbb R}$, where
\begin{equation}
\zeta_\mathrm{NLS}(Z) = \left(\frac{2a_2}{a_3}\right)^{\!\!1/2} \sech \left(\left(\frac{a_2}{a_1}\right)^{\!\!1/2}Z\right) \label{Explicit NLS}
\end{equation}
and $a_1$, $a_2$, $a_3$ are positive constants which depend upon $\omega$.
\end{theorem}

\begin{figure}[h]
\centering
\includegraphics[scale=0.5]{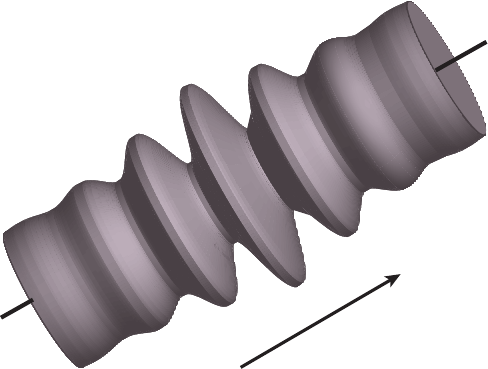}\hspace{1in}\includegraphics[scale=0.5]{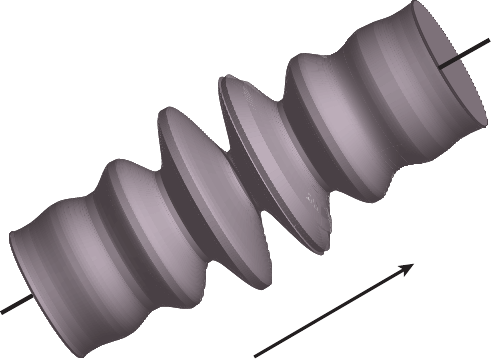}
\caption{Nonlinear Schr\"{o}dinger solitary waves of elevation (left) and of depression (right) depending on the sign in equation \eqref{eta from zeta}}
\end{figure}

Axisymmetric solitary waves have also been investigated using model equations by Bashtovoi, Rex \& Foiguel \cite{BashtovoiRexFoiguel83} and
Rannacher \& Engel \cite{RannacherEngel06}, experimentally by Bourdin, Bacri \& Falcon \cite{BourdinBacriFalcon10} and numerically by Blyth \& Parau \cite{BlythParau14}, Guyenne \& Parau \cite{GuyenneParau16},
Doak \& Vanden-Broeck \cite{DoakVandenBroeck19} and Xu \& Wang \cite{XuWang25}. Furthermore, using spatial dynamical-systems methods Groves \& Nilsson
\cite{GrovesNilsson18} have given a rigorous existence theory for multiple types of solitary waves (including those in Theorems \ref{KdV thm}
and \ref{NLS thm}).

\subsection{Weakly nonlinear theory}
\begin{figure}[h]
\centering

\includegraphics[scale=0.6]{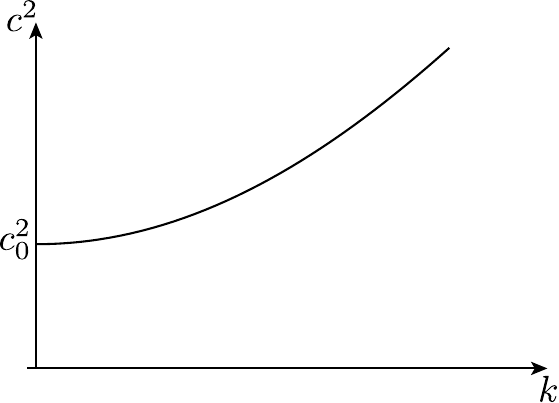}
\hspace{1.5cm}\includegraphics[scale=0.6]{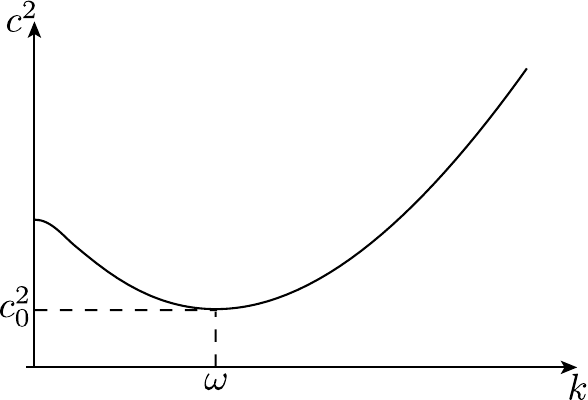}
{\it
\caption{Dispersion relation in the cases $1<\gamma\leq 9$ (left) and $\gamma>9$ (right); the minimum value of $c^2$ is denoted by $c_0^2$
\label{Dispersion relation}}}
\end{figure}

It is instructive to present a heuristic argument as a motivation for Theorems \ref{KdV thm} and \ref{NLS thm}, beginning with the linearised problem. Linearising equation \eqref{GZCS} yields
\begin{equation}
(\gamma-1)\eta - \eta_{zz} -c^2K_0 \eta=0, \label{Linear GZCS}
\end{equation}
where $K_0=f(D)$ and
$$f(|k|)=\frac{|k|I_0(|k|)}{I_1(|k|)}.$$
here we used the notation
$$h(D) \xi = \FF[h(k)\hat{\xi}], \qquad \hat{\xi}=\FF[\xi],$$
for the Fourier multiplier defined by $h$, where $\FF$ is the one-dimensional Fourier transform defined by
$$\FF[\xi](k) = \frac{1}{\sqrt{2\pi}}\int_{\mathbb{R}} \xi(z) \e^{-\i k z}\dz$$
and $D=-\i \partial_z$. Seeking solutions of \eqref{Linear GZCS} of the form $\eta(z)=\cos(kz)$ (`sinusoidal wave trains'), we obtain the dispersion
relation
$$c^2 = \frac{\gamma-1+k^2}{f(k)},$$
which describes the relation between the wave number $k \geq 0$ and the wave speed $c \geq 0$. In Appendix \ref{disprel} we show that
$c^2$ is a strictly monotone increasing function of $k$ for $1<\gamma\leq 9$, while for $\gamma>9$ it has a unique local maximum at $k=0$
and a unique global minimum at $k=\omega>0$ (the formula $\gamma=1-\omega^2 +2\omega f(\omega)/f^\prime(\omega)$ defines a bijection
between the values of $\gamma \in (9,\infty)$ and $\omega \in (0,\infty)$). In both cases we denote its global minimum by $c_0^2$, so that
$$c_0^2=\left\{\hspace{-1mm}\begin{array}{ll} c^2(0)=\frac{1}{2}(\gamma-1),\quad & 1 < \gamma \leq 9, \\[2mm]c^2(\omega)=\dfrac{2\omega}{f(\omega)}, & \gamma>9  \end{array}\right.$$
(see Figure \ref{Dispersion relation}).

Using $c$ as a bifurcation parameter, we expect branches of small-amplitude solitary waves to bifurcate
at $c=c_0$ (where the linear group and phase speeds are equal) into the region $\{c<c_0\}$ where linear periodic
wave trains are not supported (see Dias \& Kharif \cite[\S 3]{DiasKharif99}).
In the case $1<\gamma<9$, one writes $c^2=c_0^2(1-\varepsilon^2)$, 
where $\varepsilon$ is a small positive number, substitutes the Ansatz
\begin{equation}
\eta(z) = \varepsilon^2 \zeta_1(Z)+\varepsilon^4 \zeta_2(Z) + \cdots, \label{KDV Ansatz}
\end{equation}
where $Z=\varepsilon z$, into equation \eqref{GZCS}, and finds that $\zeta_1$ satisfies the
stationary Korteweg-de Vries equation
\begin{equation}
(\tfrac{1}{8}\gamma-\tfrac{9}{8})\zeta_{ZZ}+2c_0^2\zeta+2c_0^2d_0\zeta^2=0, \label{Basic KdV}
\end{equation}
which has the explicit (symmetric) solitary-wave solution $\zeta_\mathrm{KdV}$ given in Theorem \ref{KdV thm}.
In the case $\gamma>9$, one writes $c^2 =c_0^2(1-\varepsilon^2)$, uses the Ansatz
\begin{equation}
\eta(z) = \tfrac{1}{2}\varepsilon \big(\zeta_1(Z) \e^{\i \omega x}
+  \overline{\zeta_1(Z)}\e^{-\i \omega x} \big)
+\varepsilon^2 \zeta_0(Z) +\tfrac{1}{2}\varepsilon^2 \big(\zeta_2(Z) \e^{2\i \omega z}
+\overline{\zeta_2(Z)}\e^{-2\i \omega z}\big) + \cdots, \label{NLS Ansatz}
\end{equation}
where $Z=\varepsilon z$ and $\gamma=1-\omega^2 +2\omega f(\omega)/f^\prime(\omega)$,
and finds that $\zeta_1$ satisfies the stationary nonlinear Schr\"{o}dinger equation
\begin{equation}
-a_1 \zeta_{ZZ}+a_2 \zeta -a_3|\zeta|^2\zeta =0, \label{Basic NLS}
\end{equation}
which has the (symmetric) solitary-wave solutions $\pm\zeta_\mathrm{NLS}$ given in Theorem \ref{NLS thm}. Details of these
calculations are given in Appendix \ref{weakly nonlinear appendix}.

\subsection{Methodology} \label{Methodology}
\begin{figure}[h]
\centering

\includegraphics[scale=0.75]{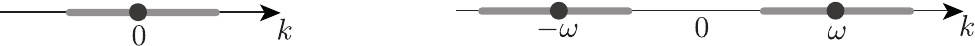}
{\it
\caption{(a) The support of $\hat{\eta}_1$ is contained in the set $S$, where $S=(-\delta,\delta)$ for $1<\gamma<9$ (left) and\linebreak
$S=(-\omega-\delta,-\omega+\delta) \cup (\omega-\delta,\omega+\delta)$ for $\gamma>9$ (right). \label{Splitting}}}
\end{figure}

In this paper we rigorously confirm the results of the weakly nonlinear theory described above.
The Ans\"{a}tze \eqref{KDV Ansatz} and \eqref{NLS Ansatz} suggest that the Fourier transform of a solitary wave 
is concentrated near the points $k=\pm \omega$ (which coincide at $k=0$ when $1<\gamma<9$). 
Indeed, writing $c^2=c_0^2(1-\varepsilon^2)$, one finds that the linearisation of \eqref{GZCS}
at $\varepsilon=0$ is
$$g(D)\eta=0,$$
where
$$
g(k):=\gamma-1+k^2-c_0^2f(k) \geq 0, \qquad k \in {\mathbb R},
$$
with equality precisely when $k=\pm \omega$ (so that $g(\omega)=g^\prime(\omega)=0$ and $g^{\prime\prime}(\omega)>0$). 
We therefore decompose $\eta$ into the sum of functions $\eta_1$ and $\eta_2$ whose Fourier transforms
$\hat{\eta}_1$ and $\hat{\eta}_2$ are supported in the region\linebreak
$S=(-\omega-\delta,-\omega+\delta) \cup (\omega-\delta,\omega+\delta)$ (with $\delta \in (0,\frac{\omega}{3})$)
and its complement (see Figure \ref{Splitting}), so that $\eta_1 = \chi(D)\eta$,
$\eta_2 = (1-\chi(D))\eta$, where $\chi$ is the characteristic function of the set $S$ (note that $S=(-\delta,\delta)$ if $\omega=0$).
Decomposing \eqref{GZCS} into
\begin{align*}
\chi(D)\left(\KK(\eta_1+\eta_2)-c_0^2(1-\varepsilon^2)\LL(\eta_1+\eta_2)\right)&=0, \\
(1-\chi(D))\left(\KK(\eta_1+\eta_2)-c_0^2(1-\varepsilon^2)\LL(\eta_1+\eta_2)\right)&=0,
\end{align*}
one finds that the second equation can be solved for $\eta_2$ as a function of $\eta_1$ for sufficiently
small values of $\varepsilon>0$; substituting 
$\eta_2=\eta_2(\eta_1)$ into the first yields the reduced equation
$$\chi(D)\left(\KK(\eta_1+\eta_2(\eta_1))-c_0^2(1-\varepsilon^2)\LL(\eta_1+\eta_2(\eta_1))\right)=0$$
for $\eta_1$ (see Section \ref{Reduction}).

Finally, the scaling
\begin{equation}
\eta_1(z) = \varepsilon^2\zeta(Z), \qquad Z=\varepsilon z, \label{KdV scaling - intro}
\end{equation}
transforms the reduced equation into
\begin{equation}
\varepsilon^{-2}g(\varepsilon D)\zeta + 2c_0^2\zeta + 2c_0^2d_0\chi_0(\varepsilon D)\zeta^2 + O(\varepsilon^{1/2})=0
\label{KdV zeta eqn - intro}
\end{equation}
for $1< \gamma < 9$, while the scaling
\begin{equation}
\eta_1(z) = \tfrac{1}{2}\varepsilon \zeta(Z) \e^{\i \omega z} + \tfrac{1}{2}\varepsilon \overline{\zeta(Z)}\e^{-\i \omega z}, \qquad Z=\varepsilon z,
\label{NLS scaling - intro}
\end{equation}
transforms the reduced equation into
\begin{equation}
\varepsilon^{-2}g(\omega+\varepsilon D)\zeta + a_2\zeta - a_3\chi_0(\varepsilon D)(|\zeta|^2\zeta) + O(\varepsilon^{1/2})=0
\label{NLS zeta eqn - intro}
\end{equation}
for $\gamma>9$; here $\chi_0$ is the characteristic function of the set $(-\delta,\delta)$, the symbol $D$ now means
$-\i \partial_Z$ and precise
estimates for the remainder terms are given in Section \ref{Reduction}.
Equations \eqref{KdV zeta eqn - intro} and \eqref{NLS zeta eqn - intro} are \emph{full dispersion} versions of (perturbed) stationary Korteweg-de Vries and nonlinear
Schr\"{o}dinger equations since they retain the linear part of the original equation \eqref{GZCS};
the fully reduced model equations \eqref{Basic KdV} and \eqref{Basic NLS}
are recovered from them in the formal limit $\varepsilon \rightarrow 0$.

The functions $\zeta_\mathrm{KdV}$ and $\pm\zeta_\mathrm{NLS}$ are nondegenerate solutions of \eqref{Basic KdV} and \eqref{Basic NLS} in the sense that the only bounded
solutions of their linearisations at $\zeta_\mathrm{KdV}$ and
$\pm\zeta_\mathrm{NLS}$ are respectively $\zeta_{\mathrm{KdV},Z}$ and $\pm\zeta_{\mathrm{NLS},Z}$, $\pm\mathrm{i}\zeta_\mathrm{NLS}$. Equation \eqref{GZCS}
is invariant under the reflection $\eta(z) \mapsto \eta(-z)$, and the reduction procedure preserves this property: the reduced
equation for $\eta_1$ is invariant under the reflection $\eta_1(z) \mapsto \eta_1(-z)$, so that \eqref{Basic KdV}
and \eqref{Basic NLS} are invariant under respectively $\zeta(Z) \mapsto \zeta(-Z)$ and $\zeta(Z) \mapsto \overline{\zeta(-Z)}$.
Restricting to spaces of symmetric functions thus eliminates the antisymmetric solutions
$\zeta_{\mathrm{KdV},Z}$ and $\pm\zeta_{\mathrm{NLS},Z}$, $\pm\mathrm{i}\zeta_\mathrm{NLS}$
of the linearised equations, and in Section \ref{sec:existence} solutions to \eqref{KdV zeta eqn - intro}
and \eqref{NLS zeta eqn - intro} are constructed as perturbations of $\zeta_\mathrm{KdV}$ and $\pm\zeta_\mathrm{NLS}$ by formulating them as fixed-point
equations and using an appropriate version of the implicit-function theorem.

This method has been used for the classical water-wave problem by Groves \cite{Groves21}, and since many of the details in the derivation and solution of our reduced equations
are similar to those in that reference we keep Sections \ref{Reduction} and \ref{sec:existence} concise. We begin our analysis by
showing that the functionals $\KK$ and $\LL$ in equation \eqref{GZCS} depend analytically upon $\eta$ in a suitable sense (see Buffoni \& Toland \cite{BuffoniToland} for
a treatise on analytic functions in Banach spaces), which of course 
entails rigorously defining the operator $K$ given by \eqref{Defn of K} and demonstrating its analyticity. This step, the details of which are given
in Section \ref{Anal}, differs significantly from the corresponding step in reference \cite{Groves21}; in particular it is necessary to study an 
axisymmetric boundary-value problem using novel function spaces and carefully estimate a Green's function defined
in terms of modified Bessel functions.

\subsection{Function spaces}

In addition to the familiar Sobolev spaces
$$H^s({\mathbb R})=\left\{\eta \in \mathscr{S}^\prime({\mathbb R}) \st \|\eta\|_s^2:=\int_{\mathbb R} (1+k^2)^s |\hat{\eta}(k)|^2 \dk < \infty\right\},\qquad s \geq 0$$
we use the variants
$$H_\varepsilon^s({\mathbb R}) = \chi_0(\varepsilon D)H^s({\mathbb R}), \qquad s \geq 0$$
and
$$
\ZZ=\left\{\eta \in \mathscr{S}^\prime({\mathbb R}) \st \|\eta\|_\ZZ := \|\hat{\eta}_1\|_{L^1({\mathbb R})} + \|\eta_2\|_2 < \infty\right\},$$
where
$$\eta_1=\chi(D) \eta, \qquad \eta_2=(1-\chi(D)) \eta$$
(see Section \ref{Methodology} above). Note in particular the estimate
$$
\|\eta_1\|_{j,\infty} \leq \|k^j\hat{\eta}(k)\|_{L^1({\mathbb R})} \lesssim \|\hat{\eta}_1\|_{L^1({\mathbb R})}, 
$$
which holds because $\hat{\eta}_1$ has compact support, and implies in particular that
\begin{equation}
\|\eta\|_{1,\infty} \leq \|\eta_1\|_{1,\infty} + \|\eta_2\|_{1,\infty} \lesssim \|\hat{\eta}_1\|_{L^1({\mathbb R})} + \|\eta_2\|_2 = \|\eta\|_\ZZ.
\label{How to estimate eta 1 infty}
\end{equation}
Our analyticity result for the operator $K$ defined by equation \eqref{Defn of K} is given in terms of the space $\ZZ$.

\begin{lemma}
The mapping $K: \ZZ \rightarrow \LL(H^{3/2}({\mathbb R}),H^{1/2}({\mathbb R}))$ is analytic at the origin.
\end{lemma}

This lemma is proved in Section \ref{Anal}, where we work with the equivalent definition
\begin{equation}
K(\eta)\xi = -(\tilde{\phi}|_{r=1+\eta})_z, \label{Alt defn of K}
\end{equation}
where $\tilde{\phi}$ is the axisymmetric solution of the Neumann boundary-value problem
\begin{align*}
\Delta \tilde{\phi} &=0, \quad\qquad 0<r<1+\eta, \\
(1+\eta)(1+\eta_z^2)^{1/2}\frac{\partial \tilde{\phi}}{\partial n}& = \xi_z, \quad\qquad r=1+\eta
\end{align*}
(which is unique up to additive constants). To solve this boundary-value problem it is obviously necessary to
study axisymmetric functions in the ferrofluid domain $\{0<r<1+\eta\}$. For this purpose we use
the radial function spaces introduced by Groves \& Hill \cite{GrovesHill24} for functions defined on the reference domain
$\{0<r<1\}$ (onto which the ferrofluid domain is mapped for our analysis). Let $\tilde{f}_m: B_1(\mathbf{0}) \times {\mathbb R} \rightarrow {\mathbb C}$ be a function with the property that
\begin{equation}
\tilde{f}_m(r\cos\theta,r\sin\theta,z)=\e^{\i m \theta}f_m(r,z), \qquad r \in [0,1),\ \theta \in {\mathbb T}^1,\ z \in {\mathbb R}, \label{eq:mode k defn}
\end{equation}
for some $m \in {\mathbb Z}$ and some $f_m:[0,1) \times {\mathbb R}\rightarrow {\mathbb C}$ with $f_m(0,z)=0$ for $m \neq 0$.
We refer to such functions as \emph{mode $m$ functions}, such that axisymmetric functions are mode $0$ functions.

\begin{remarks} \label{rem:radialcoeffs}
$ $
\begin{itemize}
\item[(i)]
The \underline{radial coefficient} $f_0(r,z)$ of a mode $0$ function $\tilde{f}_0(x,y,z)$ obviously satisfies
$f_0(0,z)=\tilde{f}_0(\mathbf{0},z)$.
The same is true for $m \neq 0$ since $f_m(0,z)=0$ implies that $\tilde{f}_m(\mathbf{0},z)=0$.
\item[(ii)]
The \underline{radial coefficient} $f_m(r,z)$ of the mode $m$ function $\tilde{f}_m(x,y,z)$ is also the radial coefficient of the
mode $-m$ function $\tilde{f}_m(x,-y,z)$.
\end{itemize}
\end{remarks}

\begin{figure}[h]
    \centering
    \includegraphics[width=0.5\textwidth]{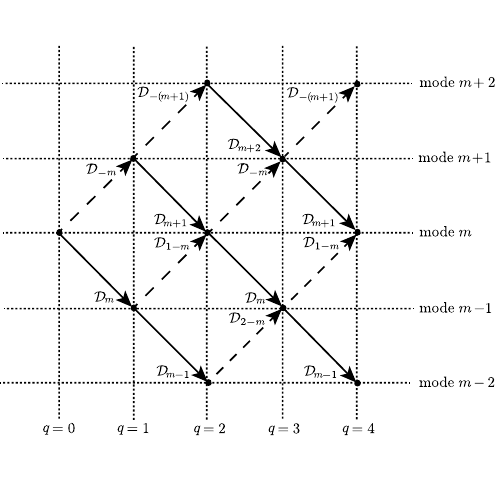}
    \caption{Actions of the Bessel operators. The mode $m+j$ function in column $q$ is $\e^{\i(m+j)\theta}{\mathcal D}_{m_q} {\mathcal D}_{m_{q-1}} \ldots {\mathcal D}_{m_1} f_m$,
    where the indices $\{m_i\}_{i=1}^{q}$ satisfy $m_i=\pm (m_{i-1} - 1)$, with $m_1=\pm m$, and consist of $\frac{1}{2}(q-j)$ positive and $\frac{1}{2}(q+j)$ non-positive terms.}
    \label{fig:Fkspace}
\end{figure}

It is convenient to study mode $m$ functions using the Wirtinger-type complex differential operators
$$
\partial_{\tau} := \tfrac{1}{\sqrt{2}}(\partial_{x} - \mathrm{i}\partial_{y}), \qquad \partial_{\bar{\tau}} := \tfrac{1}{\sqrt{2}}(\partial_{x} + \mathrm{i}\partial_{y})
$$
in place of the Cartesian differential operators $\partial_x$, $\partial_y$.
Let $\tilde{f}_m: B_1(\mathbf{0}) \times {\mathbb R} \to {\mathbb C}$ be a mode $m$ function with radial coefficient
$f_m:[0,1) \times {\mathbb R} \rightarrow {\mathbb C}$.
It follows that
$$
       \partial_{\tau}\tilde{f}_m= \mathrm{e}^{\mathrm{i}(m-1)\theta}\tfrac{1}{\sqrt{2}}\mathcal{D}_{m}f_m, \qquad \partial_{\bar{\tau}}\tilde{f}_m = \mathrm{e}^{\mathrm{i}(m+1)\theta}\tfrac{1}{\sqrt{2}}\mathcal{D}_{-m}f_m,
$$
where $ \mathcal{D}_{j}$ is the Bessel operator
$$
    \mathcal{D}_{j} := r^{-j} \frac{\mathrm{d}}{\mathrm{d}r} r^{j} = \frac{\mathrm{d}}{\mathrm{d}r} + \frac{j}{r}.
$$
According to this calculation the operators
$\partial_\tau$ and $\partial_{\bar\tau}$ map a mode $m$ function with radial coefficient $f_m$ to
a mode $m-1$ function with radial coefficient ${\mathcal D}_m f_m$ and a mode $m+1$ function with radial coefficient ${\mathcal D}_{-m} f_m$  respectively.
Correspondingly, one finds that ${\mathcal D}_m$ and ${\mathcal D}_{-m}$ map a mode $m$ radial coefficient to a mode $m-1$ and a mode $m+1$ radial coefficient respectively,
as illustrated diagrammatically in Figure \ref{fig:Fkspace} (which commutes).
Note that it is actually not necessary to distinguish between mode $m$ and mode $-m$
radial coefficients since the radial coefficient of the mode $m$ function $\tilde{f}_m(x,y)$ is also the radial coefficient of the mode $-m$ function $\tilde{f}_m(x,-y)$
(which explains the apparent ambiguity in this interpretation of ${\mathcal D}_0$.)

We denote the (closed) subspace of the standard Sobolev space
$$
H^q(B_1(\mathbf{0}) \times {\mathbb R}; {\mathbb C})=
\left\{\tilde{f}:B_1(\mathbf{0}) \times {\mathbb R}\rightarrow {\mathbb C}\st \| \tilde{f}\|_{H^q}^2 := \sum_{p=0}^{q}\sum_{n=0}^p \sum_{i=0}^n \begin{pmatrix} n \\ i \end{pmatrix}\
\| \partial_{\bar{\tau}}^{n-i} \partial_\tau^i \partial_z^{p-n}\tilde{f}\|_{L^2}^2<\infty\right\}
$$
consisting of mode $m$ functions by $\tilde{H}_{(m)}^q(B_1(\mathbf{0}) \times {\mathbb R};\mathbb{C})$.
Observe that a mode $m$ function
$\tilde{f}_m$ belongs to\linebreak $L^2(B_1(\mathbf{0}) \times {\mathbb R};{\mathbb C})$ if and only if its radial coefficient $f_m$ belongs to
$$L_1^2((0,1)\times{\mathbb R};{\mathbb C}) = \left\{f: [0,1)\times{\mathbb R} \rightarrow {\mathbb C}\ \big|\  \|f\|_{L_1^2}^2:=2\pi\int_{{\mathbb R}} \int_0^1 |f(r,z)|^2 r \dr \dz< \infty\right\},$$
and that
$$
\partial_{\bar{\tau}}^{n-i}\partial^{i}_{\tau}\tilde{f}_m = \mathrm{e}^{\mathrm{i}(m+n-2i)\theta}\,2^{-\frac{n}{2}}\mathcal{D}_{-m+i}^{n-i}\mathcal{D}^{i}_{m}f_m,
\qquad i=0,\ldots n,
$$
where
$$
    \mathcal{D}^{i}_{j} := r^{-j+i} \left(\frac{1}{r}\frac{\mathrm{d}}{\mathrm{d}r}\right)^{i} r^{j} = \mathcal{D}_{j-(i-1)}\mathcal{D}_{j-(i-2)}\dots\mathcal{D}_{j-1}\mathcal{D}_{j}.
$$
One indeed finds that $\tilde{f}_m$ belongs to $\tilde{H}_{(m)}^q(B_1(\mathbf{0}) \times {\mathbb R};{\mathbb C})$ if and only if
its radial coefficient $f_m$ belongs to 
$$
H_{(m)}^q((0,1) \times {\mathbb R};{\mathbb C})\!=\!\left\{\!f_m: [0,1) \times {\mathbb R} \rightarrow {\mathbb C}\ \big|\ \|f_m\|_{H^{q}_{(m)}}^2\!\!\!:= \sum_{p=0}^{q}\sum_{n=0}^p 2^{-n} \sum_{i=0}^n \begin{pmatrix} n \\ i \end{pmatrix} 
\|{\mathcal D}_{-m+i}^{n-i}{\mathcal D}_m^i \partial_z^{p-n}f_m\|_{L_1^2}^2<\infty\!\right\}\!,
$$
and that the mapping 
$f_m \mapsto \tilde{f}_m$ is an isometric isomorphism (see Groves \& Hill \cite[\S3]{GrovesHill24} for a more precise statement and a discussion
of the properties of these function spaces).

\section{Analyticity} \label{Anal}

\subsection{The operator $K$} \label{Everything is analytic}

In this section we study the operator $K$ given by \eqref{Alt defn of K}.
Denoting the radial coefficient of $\tilde{\phi}$ by $\phi$, such that
$$\tilde{\phi}(x,y,z)=\phi(r,z),$$
we can equivalently define
$$K(\eta)\xi = -(\phi|_{r=1+\eta})_z,$$
where $\phi$ is the solution of the boundary-value problem
\begin{align}
{\mathcal D}_1 {\mathcal D}_0 \phi + \phi_{zz}  &=0, \quad\qquad 0<r<1+\eta, \label{BVP for phi 1} \\
(1+\eta)( {\mathcal D}_0\phi -\eta_z\phi_z)&=\xi_z, \quad\qquad\!\! r=1+\eta \label{BVP for phi 2}
\end{align}
(which is unique up to additive constants).

The `flattening' transformation
$$r^\prime=\frac{r}{1+\eta}, \qquad u(r^\prime,z)=\phi(r,z)$$
transforms $S_1$ into the fixed strip $\Sigma=(0,1) \times {\mathbb R}$ and the boundary-value problem \eqref{BVP for phi 1},
\eqref{BVP for phi 2} into
\begin{alignat}{2}
{\mathcal D}_1 {\mathcal D}_0 u + u_{zz} & = {\mathcal D}_1 F_1(\eta,u) + \partial_z F_2(\eta,u), \qquad && 0<r<1, \label{Flattened BC 1} \\
{\mathcal D}_0 u & = F_1(\eta,u)+\xi_z, && r=1, \label{Flattened BC 2}
\end{alignat}
where we have dropped the primes for notational simplicity and
\begin{equation}
F_1(\eta,u)=r(1+\eta)\eta_zu_z - r^2\eta_z^2{\mathcal D}_0u, \qquad F_2(\eta,u)=r(1+\eta)\eta_z {\mathcal D}_0u - \eta(\eta+2)u_z, \label{Defn of F1, F2}
\end{equation}
so that
$$K(\eta)\xi=-u_z|_{r=1}.$$ 
This boundary-value problem can be cast as the integral equation
\begin{equation}
u = S(F_1(\eta,u),F_2(\eta,u),\xi), \label{Flattened integral eqn}
\end{equation}
where
$$S(F_1,F_2,\xi) = {\mathcal F}^{-1}\left[\int_0^1 \left(\i k G(r,\tilde{r})\hat{F}_2 - \tilde{\mathcal D}_0 G(r,\tilde{r}) \hat{F}_1 \right)\tilde{r}\dtilder - \i k G(r,1)\hat{\xi}\right]$$
with
$$G(r,\tilde{r}) = \left\{\begin{array}{ll}
-I_0(|k|r)\left(K_0(|k|\tilde{r})+\dfrac{K_1(|k|)}{I_1(|k|)} I_0(|k|\tilde{r})\right), & 0 \leq r < \tilde{r}, \\[4mm]
-I_0(|k|\tilde{r})\left(K_0(|k|r)+\dfrac{K_1(|k|)}{I_1(|k|)} I_0(|k|r)\right), & \tilde{r} < r < 1.
\end{array}\right.$$

We study \eqref{Flattened integral eqn} for $\eta \in \ZZ$, $\xi \in H^{3/2}({\mathbb R})$ and $u \in H^\star(\Sigma)$, where
$$H^\star(\Sigma):=H_{(0)}^2(\Sigma)/{\mathbb R}$$ with norm
$$\|u\|_\star^2 := \|u_z\|_{H_{(0)}^1}^2 + \|{\mathcal D}_0 u\|_{H_{(1)}^1}^2.$$
The following result is proved in Section \ref{Linear BVP} below.

\begin{theorem} \label{Estimates for S}
The solution operator $S$ satisfies
$$\|S(F_1,F_2,\xi)\|_\star \lesssim \|F_1\|_{H^1_{(1)}}+\|F_2\|_{H^1_{(0)}}+\|\xi\|_{3/2}$$
for all $F_1 \in H^1_{(1)}(\Sigma)$, $F_2 \in H^1_{(0)}(\Sigma)$ and $\xi \in H^{3/2}({\mathbb R})$.
\end{theorem}

\begin{lemma} \label{F1 and F2 are analytic}
The formulae \eqref{Defn of F1, F2} define analytic functions $F_1: \ZZ \times H^\star(\Sigma) \rightarrow H_{(1)}^1(\Sigma)$,
$F_2: \ZZ \times H^\star(\Sigma) \rightarrow H_{(0)}^1(\Sigma)$.
\end{lemma}
{\bf Proof.} 
Clearly
$$\|r\eta_zu_z\|_{L_1^2} \lesssim \|\eta_z\|_\infty \|u_z\|_{L_1^2} \lesssim \|\eta\|_\ZZ \|u\|_\star,$$
(see \eqref{How to estimate eta 1 infty}).
Using the calculations
\begin{alignat*}{2}
& {\mathcal D}_1(r \eta_z u_z) = r \eta_z {\mathcal D}_0(u_z)+2\eta_zu_z, \qquad \partial_z(r \eta_z u_z)&&=r\eta_{zz}u_z+r\eta_z u_{zz} \\
& && = r\eta_{1zz}u_z + r\eta_{2zz}u_z + r \eta_z u_{zz},
\end{alignat*}
we similarly find that
$$\|{\mathcal D}_1(r \eta_z u_z)\|_{L_1^2} \lesssim \|\eta_z\|_\infty (\|{\mathcal D}_0 u_z\|_{L_1^2} + \|u_z\|_{L_1^2}) \lesssim \|\eta\|_\ZZ  \|u\|_\star$$
and
$$\|\partial_z(r \eta_z u_z)\|_{L_1^2} \lesssim \|\eta_{1zz}\|_\infty\|u_z\|_{L_1^2} + \|\eta_z\|_\infty \|u_{zz}\|_{L_1^2} + \|\eta_{2zz}u_z\|_{L_1^2}
\lesssim \|\eta\|_\ZZ  \|u\|_\star + \|\eta_{2zz}u_z\|_{L_1^2}$$
with
\begin{align*}
\|\eta_{2zz} u_z \|_{L_1^2}^2 & \leq \|\eta_{2zz}\|_0^2 \sup_{z \in {\mathbb R}}\|r^{1/2}u_z\|^2_{L^2(0,1)} \\
& \lesssim \|\eta_{2zz}\|_0^2\|r^{1/2} u_z(r,z)\|_{H^1({\mathbb R},L^2(0,1))}^2 \\
& =\|\eta_{2zz}\|_0^2(\|u_z\|_{L_1^2}^2+\|u_{zz}\|_{L_1^2}^2) \\
&\leq \|\eta\|_\ZZ^2 \|u\|_\star^2.
\end{align*}
It follows that $(\eta,u) \mapsto r\eta_zu_z$ is an analytic mapping $\ZZ \times H^\star(\Sigma) \rightarrow H_{(1)}^1(\Sigma)$, and similar arguments show\linebreak that
$(\eta,u) \mapsto r\eta\eta_zu_z$, $(\eta,u) \mapsto r^2\eta_z^2{\mathcal D}_0u$ are analytic mappings $\ZZ \times H^\star(\Sigma) \rightarrow H_{(1)}^1(\Sigma)$,
such that\linebreak $F_1: \ZZ \times H^\star(\Sigma) \rightarrow H_{(1)}^1(\Sigma)$ is analytic. The same method shows that $F_2: \ZZ \times H^\star(\Sigma) \rightarrow H_{(0)}^1(\Sigma)$ is analytic.\qed

\begin{theorem}
For each $\xi\in H^{3/2}({\mathbb R})$ and each sufficiently small $\eta \in \ZZ$ the boundary-value problem \eqref{Flattened BC 1},
\eqref{Flattened BC 2} admits a unique solution $u \in H^\star(\Sigma)$. Furthermore, the mapping $\ZZ \mapsto \LL(H^{3/2}({\mathbb R}),H^\star(\Sigma))$
is analytic at the origin.
\end{theorem}
{\bf Proof.} Define a mapping $T: H^\star(\Sigma) \times \ZZ \times H^{3/2}({\mathbb R}) \rightarrow H^\star(\Sigma)$ by
$$T(u,\eta,\xi) = u - S(F_1(\eta,u),F_2(\eta,u),\xi),$$
such that the solutions to \eqref{Flattened integral eqn} are precisely the zeros of $T(\cdot,\eta,\xi)$. It follows from
Theorem \ref{Estimates for S} and Lemma \ref{F1 and F2 are analytic} that $T$ is analytic at the origin. Furthermore $T(0,0,0)=0$
and $\mathrm{d}_1T[0,0,0]=I$ is an isomorphism. It follows from the analytic implicit-function theorem (see Buffoni \& Toland
\cite[Theorem 4.5.4]{BuffoniToland}) that there exists open neighbourhoods $N_1 \subseteq \ZZ$, $N_2 \subseteq H^{3/2}({\mathbb R})$
and $N_3 \subseteq H^\star(\Sigma)$ of the origin and analytic function $v: N_1 \times N_2 \rightarrow N_3$ such that
$$T(v(\eta,\xi),\eta,\xi)=0;$$
furthermore $u=v(\eta,\xi)$ for all $(\eta,\xi,u) \in N_1 \times N_2 \times N_3$ with $T(u,\eta,\xi)=0$. Since $u$ is linear in $\xi$ we can
choose $N_2=H^{3/2}({\mathbb R})$.\qed
\begin{corollary} \label{K is analytic}
The mapping $K: \ZZ \rightarrow \LL(H^{3/2}({\mathbb R}),H^{1/2}({\mathbb R}))$ is analytic at the origin.
\end{corollary}
{\bf Proof.} This assertion follows from the formula $K(\eta)\xi = -u_z|_{r=1}$, the analyticity of $u: N_1 \times H^{3/2}({\mathbb R}) \rightarrow H^\star(\Sigma)$ and the facts that $\partial_z: H^\star(\Sigma) \rightarrow H_{(0)}^1(\Sigma)$ and $u \mapsto u|_{r=1}$, $H_{(0)}^1(\Sigma) \rightarrow H^{1/2}({\mathbb R})$ are continuous linear operators (see Groves \& Hill \cite[Lemma 3.24]{GrovesHill24}).\qed

According to Corollary \ref{K is analytic} we can choose $M$ sufficiently small and study the equation
$$\KK(\eta)-c_0^2(1-\varepsilon^2)\LL(\eta)=0$$
in the set
\begin{equation}
U=\{\eta \in H^2({\mathbb R}): \|\eta\|_\ZZ < M\}, \label{Defn of U}
\end{equation}
noting that $H^2({\mathbb R})$ is continuously embedded in $\ZZ$ and $U$ is an open neighbourhood of the origin in $H^2({\mathbb R})$.

\begin{corollary} \label{KK and LL are analytic}
The formulae \eqref{Defn of KK}, \eqref{Defn of LL} define analyic functions $U \rightarrow L^2({\mathbb R})$.
\end{corollary}
{\bf Proof.} This observation follows from the formulae
\begin{align*}
\KK(\eta) & = -\gamma\left(\nu\left(\frac{1}{1+\eta}\right)-\nu(1)\right)+\left(\frac{1}{1+\eta}-1\right)\!\!\left(\frac{1}{(1+\eta_z^2)^{1/2}}-1\right)+\frac{1}{1+\eta}-1 \\
& \quad\qquad\mbox{}+\frac{1}{(1+\eta_z^2)^{1/2}}-1-\left(\frac{1}{(1+\eta_z^2)^{3/2}}-1\right)\eta_{zz}-\eta_{zz}, \\
\LL(\eta) & = -\tfrac{1}{2}(K(\eta)\eta)^2-\tfrac{1}{2}K(\eta)\eta K(\eta)\eta^2 - \tfrac{1}{8}(K(\eta)\eta^2)^2+\frac{\eta_z^2}{2(1+\eta_z^2)}
+\frac{\eta_z^2}{2(1+\eta_z^2)}(K(\eta)\eta)^2 \\
& \quad\qquad\mbox{}+\frac{\eta_z^2}{8(1+\eta_z^2)}(K(\eta)\eta^2)^2-\frac{\eta_z^2}{1+\eta_z^2}K(\eta)\eta-\frac{\eta_z^2}{2(1+\eta_z^2)}K(\eta)\eta^2 \\
& \quad\qquad\mbox{}+\frac{\eta_z^2}{2(1+\eta_z^2)}K(\eta)\eta K(\eta)\eta^2+K(\eta)\eta+\tfrac{1}{2}K(\eta)\eta^2
\end{align*}
and
\begin{itemize}
\item[(i)]
Corollary \ref{K is analytic},
\item[(i)]
the fact that the functions
$$
\rho \mapsto \nu\left(\frac{1}{1+\rho}\right)-\nu(1), \qquad
\rho \mapsto \frac{1}{1+\rho}-1, \qquad
\rho \mapsto \frac{1}{(1+\rho^2)^{1/2}}-1, \qquad
\rho \mapsto \frac{\rho^2}{(1+\rho^2)^{1/2}}
$$
are analytic at the origin $H^1({\mathbb R}) \rightarrow H^1({\mathbb R})$,
\item[(ii)]
the continuity of the multiplication map $H^1({\mathbb R}) \times H^1({\mathbb R}) \rightarrow L^2({\mathbb R})$,
$H^1({\mathbb R}) \times L^2({\mathbb R}) \rightarrow L^2({\mathbb R})$
and $H^{1/2}({\mathbb R}) \times H^{1/2}({\mathbb R}) \rightarrow L^2({\mathbb R})$
(see H\"{o}rmander \cite[Theorem 8.3.1]{Hoermander}),
\item[(iii)]
the continuity of the embeddings $H^2({\mathbb R}) \subseteq H^{3/2}({\mathbb R}) \subseteq H^1({\mathbb R}) \subseteq H^{1/2}({\mathbb R}) \subseteq L^2({\mathbb R})$,
\item[(iv)]
the fact that $H^{3/2}({\mathbb R})$ is a Banach algebra.\qed
\end{itemize}

\subsection{The linear boundary-value problem} \label{Linear BVP}

In this section we prove Theorem \ref{Estimates for S} by estimating the operators
\begin{align*}
\GG_1(F)&:=\FF^{-1}\left[\int_0^1 \i k \tilde{r} G(r,\tilde{r})\hat{F}(\tilde{r})\dtilder\right], \\
\GG_2(F)&:=\FF^{-1}\left[\int_0^1 -\tilde{\mathcal D}_0 G(r,\tilde{r})\tilde{r}\hat{F}(\tilde{r})\dtilder\right], \\
\GG_3(\xi)&:=\FF^{-1}[-\i k G(r,1)\hat{\xi}]=\FF^{-1}\left[-\i\frac{I_0(|k|r)}{I_1(|k|)}\hat{\xi}\right]
\end{align*}
in Lemmata \ref{Estimate for G1}--\ref{Estimate for G3} below. For this purpose we introduce the functions
\begin{align*}
H_1(r,\tilde{r}) &= \left\{\begin{array}{ll}
-|k|I_1(|k|r)\left(K_0(|k|\tilde{r})+\dfrac{K_1(|k|)}{I_1(|k|)} I_0(|k|\tilde{r})\right), & 0 \leq r < \tilde{r}, \\[4mm]
|k|I_0(|k|\tilde{r})\left(K_1(|k|r)-\dfrac{K_1(|k|)}{I_1(|k|)} I_1(|k|r)\right), & \tilde{r} < r < 1,
\end{array}\right. \\
H_2(r,\tilde{r}) &= \left\{\begin{array}{ll}
|k|I_0(|k|r)\left(K_1(|k|\tilde{r})-\dfrac{K_1(|k|)}{I_1(|k|)} I_1(|k|\tilde{r})\right), & 0 \leq r < \tilde{r}, \\[4mm]
-|k|I_1(|k|\tilde{r})\left(K_0(|k|r)+\dfrac{K_1(|k|)}{I_1(|k|)} I_0(|k|r)\right), & \tilde{r} < r < 1,
\end{array}\right.\\
H_3(r,\tilde{r}) &= \left\{\begin{array}{ll}
|k|^2I_1(|k|r)\left(K_1(|k|\tilde{r})-\dfrac{K_1(|k|)}{I_1(|k|)} I_1(|k|\tilde{r})\right), & 0 \leq r < \tilde{r}, \\[4mm]
|k|^2I_1(|k|\tilde{r})\left(K_1(|k|r)-\dfrac{K_1(|k|)}{I_1(|k|)} I_1(|k|r)\right), & \tilde{r} < r < 1,
\end{array}\right.
\end{align*}
which are the formal derivatives $G_r(r,\tilde{r})$, $G_{\tilde{r}}(r,\tilde{r})$ and $G_{r\tilde{r}}(r,\tilde{r})$ of $G(r,\tilde{r})$
respectively, and establish the following auxiliary results.

\begin{proposition} \label{First proposition}
The function $G(r,\tilde{r})$ satisfies
$$
\int_0^1 \tilde{r}|G(r,\tilde{r})|\dtilder = \frac{1}{|k|^2}, \qquad
\int_0^1 r|G(r,\tilde{r})|\dr = \frac{1}{|k|^2}
$$
for all $k \in {\mathbb R}$.
\end{proposition}
{\bf Proof.} We find that
\begin{align*}
\int_0^1 r|G(r,\tilde{r})|\dr & = \left(K_0(|k|\tilde{r})+\frac{K_1(|k|)}{I_1(|k|)}I_0(|k|\tilde{r})\right)\int_0^{\tilde{r}} rI_0(|k|r) \dr \\
& \quad\qquad\mbox{}+I_0(|k|\tilde{r}) \int_{\tilde{r}}^1 r\left( K_0(|k|r)+\frac{K_1(|k|)}{I_1(|k|)}I_0(|k|r)\right)\dr \\
& = \frac{1}{|k|^2}
\end{align*}
and
$$\int_0^1 \tilde{r}|G(r,\tilde{r})|\dtilder = \int_0^1 r|G(r,\tilde{r})|\dr = \frac{1}{|k|^2}.\eqno{\Box}$$

\begin{proposition} \label{Second proposition}
The function $H_1(r,\tilde{r})$ satisfies
$$
\int_0^1 \tilde{r}|H_1(r,\tilde{r})|\dtilder \lesssim \frac{1}{|k|}, \qquad
\int_0^1 r|H_1(r,\tilde{r})|\dr \lesssim \frac{1}{|k|}
$$
for all $k \in {\mathbb R}$.
\end{proposition}
{\bf Proof.} We note that
$$2|k|r I_1(|k|r)K_1(|k|r) \rightarrow \left\{\begin{array}{ll}0 \quad \mbox{as $|k|r \rightarrow 0$}, \\[2mm]
1 \quad \mbox{as $|k|r \rightarrow \infty$}\end{array}\right.$$
and
$$2|k|rI_1(|k|r)^2\frac{K_1(|k|)}{I_1(|k|)}
\leq 2|k| I_1(|k|) K_1(|k|)
\rightarrow \left\{\begin{array}{ll}0 \quad \mbox{as $|k| \rightarrow 0$}, \\[2mm]
1 \quad \mbox{as $|k| \rightarrow \infty$},\end{array}\right.$$
so that these quantities are bounded over $r \in [0,1]$ and $k \in {\mathbb R}$. We therefore find that
\begin{align*}
|k|\int_0^1 \tilde{r}|H_1(r,\tilde{r})\dtilder
& = |k|\left(K_1(|k|r)-\frac{K_1(|k|)}{I_1(|k|)}I_1(|k|r)\right)\int_0^r |k| \tilde{r}I_0(|k|\tilde{r}) \dtilder \\
& \quad\qquad\mbox{}+|k|I_1(|k|r) \int_r^1 |k|\tilde{r}\left( K_0(|k|\tilde{r})+\frac{K_1(|k|)}{I_1(|k|)}I_0(|k|\tilde{r})\right)\dtilder\\
& = -2|k|rI_1(|k|r)^2\frac{K_1(|k|)}{I_1(|k|)} + 2|k|rI_1(|k|r)K_1(|k|r) \\
& \lesssim 1.
\end{align*}
Next we record the estimates
\begin{align*}
0 & \leq \tfrac{\pi}{2}|k|\tilde{r}I_1(|k|\tilde{r})K_0(|k|\tilde{r})L_0(|k|\tilde{r}) - \tfrac{\pi}{2}|k|\tilde{r}I_0(|k|\tilde{r})K_1(|k|\tilde{r})L_0(|k|\tilde{r}) \\
& \quad\qquad\mbox{}-\pi |k|\tilde{r}I_0(|k|\tilde{r})K_0(|k|\tilde{r})L_1(|k|\tilde{r})+\tfrac{\pi}{2}I_0(|k|\tilde{r}) \\
& \rightarrow \left\{\begin{array}{ll}\frac{\pi}{2} \quad \mbox{as $|k|\tilde{r} \rightarrow 0$}, \\[2mm]
1 \quad \mbox{as $|k|\tilde{r} \rightarrow \infty$},\end{array}\right.\\
\\
0 & \leq \tfrac{\pi}{2}I_0(|k|\tilde{r}) \left(1-\frac{L_1(|k|)}{I_1(|k|)}\right) 
 \leq \tfrac{\pi}{2}I_0(|k|) \left(1-\frac{L_1(|k|)}{I_1(|k|)}\right) 
\rightarrow \left\{\begin{array}{ll}\frac{\pi}{2} \quad \mbox{as $|k| \rightarrow 0$}, \\[2mm]
1 \quad \mbox{as $|k| \rightarrow \infty$},\end{array}\right.
\end{align*}
where $L_\nu$ is the modified Struve function of the first kind and order $\nu$. Using the fact that\linebreak $h: s \mapsto \pi s(I_1(s)L_0(s)-I_0(s)L_1(s))$
is increasing (since $h^\prime(s) = 2sI_1(s) > 0$ for $s>0$ with $h^\prime(0)=0$), we furthermore find that
\begin{align*}
0 & \leq \pi|k|\tilde{r}I_0(|k|\tilde{r})\Big(I_1(|k|\tilde{r})L_0(|k|\tilde{r})-I_0(|k|\tilde{r})L_1(|k|\tilde{r})\Big)\frac{K_1(|k|\tilde{r})}{I_1(|k|\tilde{r})} \\
& \leq \pi|k|I_0(|k|)\Big(I_1(|k|)L_0(|k|)-I_0(|k|)L_1(|k|)\Big)\frac{K_1(|k|)}{I_1(|k|)} \\
& \rightarrow \left\{\begin{array}{ll}0 \quad \mbox{as $|k| \rightarrow 0$}, \\[2mm]
1 \quad \mbox{as $|k| \rightarrow \infty$}.\end{array}\right.\\
\end{align*}
Using these estimates we conclude that
\begin{align*}
|k|\int_0^1 r |H_1(r,\tilde{r})|\dr
& = |k|\left(K_0(|k|\tilde{r})+\frac{K_1(|k|)}{I_1(|k|)}I_0(|k|\tilde{r})\right)\int_0^{\tilde{r}} |k| rI_1(|k|r) \dtilder \\
& \quad\qquad\mbox{}+|k|I_0(|k|\tilde{r}) \int_{\tilde{r}}^1 |k|r\left( K_1(|k|r)-\frac{K_1(|k|)}{I_1(|k|)}I_1(|k|r)\right)\dtilder\\
& = \tfrac{\pi}{2}|k|\tilde{r}I_1(|k|\tilde{r})K_0(|k|\tilde{r})L_0(|k|\tilde{r}) - \tfrac{\pi}{2}|k|\tilde{r}I_0(|k|\tilde{r})K_1(|k|\tilde{r})L_0(|k|\tilde{r}) \\
& \quad\qquad\mbox{}-\pi |k|\tilde{r}I_0(|k|\tilde{r})K_0(|k|\tilde{r})L_1(|k|\tilde{r})+\tfrac{\pi}{2}I_0(|k|\tilde{r}) \frac{L_1(|k|)}{I_1(|k|)}\\
& \quad\qquad\mbox{}+\pi|k|\tilde{r}I_0(|k|\tilde{r})\Big(I_1(|k|\tilde{r})L_0(|k|\tilde{r})-I_0(|k|\tilde{r})L_1(|k|\tilde{r})\Big)\frac{K_1(|k|)}{I_1(|k|)} \\
& \lesssim 1.\alignqed
\end{align*}

\begin{corollary} \label{Third proposition}
The function $H_2(r,\tilde{r})$ satisfies
$$
\int_0^1 \tilde{r}|H_2(r,\tilde{r})|\dtilder \lesssim \frac{1}{|k|}, \qquad
\int_0^1 r|H_2(r,\tilde{r})|\dr \lesssim \frac{1}{|k|}
$$
for all $k \in {\mathbb R}$.
\end{corollary}

\begin{proposition} \label{Fourth proposition}
The function $H_3(r,\tilde{r})$ satisfies
$$
\int_0^1 \tilde{r}|H_3(r,\tilde{r})|\dtilder \lesssim 1, \qquad
\int_0^1 r|H_3(r,\tilde{r})|\dr \lesssim 1
$$
for all $k \in {\mathbb R}$.
\end{proposition}
{\bf Proof.} Using the estimates
\begin{align*}
0 & \leq \tfrac{\pi}{2}I_1(|k|r)\left(1-\frac{L_1(|k|)}{I_1(|k|)}\right) 
 \leq \tfrac{\pi}{2}I_1(|k|)\left(1-\frac{L_1(|k|)}{I_1(|k|)}\right) 
\rightarrow \left\{\begin{array}{ll}0 \quad \mbox{as $|k| \rightarrow 0$}, \\[2mm]
1 \quad \mbox{as $|k| \rightarrow \infty$},\end{array}\right.\\
\\
0 & \leq \tfrac{\pi}{2}\big(I_1(|k|r)-L_1(|k|r)\big) 
\rightarrow \left\{\begin{array}{ll}0 \quad \mbox{as $|k| \rightarrow 0$}, \\[2mm]
1 \quad \mbox{as $|k| \rightarrow \infty$},\end{array}\right.
\end{align*}
we find that
\begin{align*}
\int_0^1 \tilde{r}|H_3(r,\tilde{r})|\dtilder & = |k|\left(K_1(|k|r)-\frac{K_1(|k|)}{I_1(|k|)}I_1(|k|r)\right)\int_0^r |k|\tilde{r}I_1(|k|\tilde{r}) \dtilder \\
& \quad\qquad\mbox{}+|k|I_1(|k|r) \int_r^1 |k|\tilde{r}\left( K_1(|k|\tilde{r})-\frac{K_1(|k|)}{I_1(|k|)}I_1(|k|\tilde{r})\right)\dtilder \\
& = \tfrac{\pi}{2}I_1(|k|r)\frac{L_1(|k|)}{I_1(|k|)}-\tfrac{\pi}{2}L_1(|k|r) \\
& \lesssim 1
\end{align*}
and
$$\int_0^1 r|H_3(r,\tilde{r})|\dr = \int_0^1 \tilde{r}|H_3(r,\tilde{r})|\dtilder \lesssim 1.\eqno{\Box}$$

\begin{lemma} \label{Estimate for G1}
The estimate
$$\|\GG_1(F)\|_\star \lesssim \|F\|_{H_{(0)}^1}$$
holds for all $F \in H_{(0)}^1(\Sigma)$.
\end{lemma}
{\bf Proof.} It follows fom Proposition \ref{First proposition} that
\begin{align*}
\|\partial_z \GG_1(F)\|_{L_1^2}^2 & \lesssim \int_{-\infty}^\infty \int_0^1 r\left|\int_0^1 -|k|^2 \tilde{r}G(r,\tilde{r})\hat{F}(\tilde{r})\dtilder\right|^2\dr\dk \\
& \lesssim \int_{-\infty}^\infty \int_0^1 r \int_0^1 |k|^2 \tilde{r}|G(r,\tilde{r})|\dtilder  \int_0^1 |k|^2 \tilde{r}|G(r,\tilde{r})||\hat{F}(\tilde{r})|^2\dtilder \dr\dk \\
&\lesssim \int_{-\infty}^\infty \int_0^1 \left(\int_0^1 |k|^2 r |G(r,\tilde{r})|\dr \right) \tilde{r}|\hat{F}(\tilde{r})|^2\dtilder\dk \\
& \lesssim \|F\|_{L_1^2}^2
\end{align*}
and hence
$$\|\partial_z^2 \GG_1(F)\|_{L_1^2} = \| \partial_z \GG_1(F_z)\|_{L_1^2} \lesssim \|F_z\|_{L_1^2}.$$

Next we note that
$${\mathcal D}_0 \GG_1(F) = \FF^{-1}\left[\int_0^1 \i k \tilde{r}H_1(r,\tilde{r})\hat{F}(\tilde{r})\dtilder\right]$$
and using Proposition \ref{Second proposition} that
\begin{align*}
\|{\mathcal D}_0 \GG_1(F)\|_{L_1^2}^2 & \lesssim \int_{-\infty}^\infty \int_0^1 r\left|\int_0^1 \i k \tilde{r}H_1(r,\tilde{r})\hat{F}(\tilde{r})\dtilder\right|^2\dr\dk \\
& \lesssim \int_{-\infty}^\infty \int_0^1 r \int_0^1 |k| \tilde{r}|H_1(r,\tilde{r})|\dtilder  \int_0^1 |k| \tilde{r}|H_1(r,\tilde{r})||\hat{F}(\tilde{r})|^2\dtilder \dr\dk \\
& \lesssim \int_{-\infty}^\infty \int_0^1\left(\int_0^1 |k| r |H_1(r,\tilde{r})|\dr\right) \tilde{r}|\hat{F}(\tilde{r})|^2\dtilder \dk\\
& \lesssim \|F\|_{L_1^2}^2
\end{align*}
and hence
$$\|\partial_z {\mathcal D}_0\GG_1(F)\|_{L_1^2} = \| {\mathcal D}_0 \GG_1(F_z)\|_{L_1^2} \lesssim \|F_z\|_{L_1^2}.$$

Finally, using the identity
$${\mathcal D}_1{\mathcal D}_0 \GG_1(F) + \partial_z^2 \GG_1(F) = F_z,$$
which follows from the definition of $\GG_1$, we find that
$$\|{\mathcal D}_1{\mathcal D}_0 \GG_1(F)\|_{L_1^2} \leq \|\partial_z^2 \GG_1(F)\|_{L_1^2} + \|F_z\|_{L_1^2} \lesssim \|F_z\|_{L_1^2}.\eqno{\Box}.$$

\begin{lemma} \label{Estimate for G2}
The estimate
$$\|\GG_2(F)\|_\star \lesssim \|F\|_{H_{(1)}^1}$$
holds for all $F \in H_{(1)}^1(\Sigma)$.
\end{lemma}
{\bf Proof.} Since
$$\partial_z \GG_2(F) = \FF^{-1}\left[\int_0^1 \i k \tilde{r}H_2(r,\tilde{r})\hat{F}(\tilde{r})\dtilder\right]$$
it follows from Corollary \ref{Third proposition} that
\begin{align*}
\|\partial_z \GG_2(F)\|_{L_1^2}^2 & \lesssim \int_{-\infty}^\infty \int_0^1 r\left|\int_0^1 \i k \tilde{r}H_2(r,\tilde{r})\hat{F}(\tilde{r})\dtilder\right|^2\dr\dk \\
& \lesssim \int_{-\infty}^\infty \int_0^1 r |k|\int_0^1  \tilde{r}|H_2(r,\tilde{r})|\dtilder  \int_0^1 |k| \tilde{r}|H_2(r,\tilde{r})||\hat{F}(\tilde{r})|^2\dtilder \dr\dk \\
& \lesssim \int_{-\infty}^\infty \int_0^1\left( \int_0^1 |k| r |H_2(r,\tilde{r})|\dr \right) \tilde{r}|\hat{F}(\tilde{r})|^2\dtilder \dk \\
& \lesssim \|F\|_{L_1^2}^2
\end{align*}
and hence
$$\|\partial_z^2 \GG_2(F)\|_{L_1^2} = \| \partial_z \GG_2(F_z)\|_{L_1^2} \lesssim \|F_z\|_{L_1^2}.$$
{Next we note that
$${\mathcal D}_0 \GG_2(F) = \FF^{-1}\left[\int_0^1 \tilde{r}H_3(r,\tilde{r})\hat{F}(\tilde{r})\dtilder\right] - F,$$
so that
\begin{align*}
\|\mathcal{D}_0\GG_2(F)\|_{L_1^2}^2 & \lesssim \int_{-\infty}^\infty \int_0^1 r\left|\int_0^1 \tilde{r}H_3(r,\tilde{r})\hat{F}(\tilde{r})\dtilder\right|^2\dr\dk + \|F\|_{L_1^2}\\
& \lesssim \int_{-\infty}^\infty \int_0^1 r \int_0^1  \tilde{r}|H_3(r,\tilde{r})|\dtilder  \int_0^1 \tilde{r}|H_3(r,\tilde{r})||\hat{F}(\tilde{r})|^2\dtilder \dr\dk + \|F\|_{L_1^2}\\
& \lesssim \int_{-\infty}^\infty \int_0^1\left( \int_0^1  r |H_3(r,\tilde{r})|\dr \right)\tilde{r}|\hat{F}(\tilde{r})|^2\dtilder \dk\\
& \lesssim \|F\|_{L_1^2}^2,
\end{align*}
where we have used Proposition \ref{Fourth proposition}, and hence
$$\|\partial_z {\mathcal D}_0\GG_2(F)\|_{L_1^2} = \| {\mathcal D}_0 \GG_2(F_z)\|_{L_1^2} \lesssim \|F_z\|_{L_1^2}.$$

Finally, using the identity
$${\mathcal D}_1{\mathcal D}_0 \GG_2(F) + \partial_z^2 \GG_2(F) = {\mathcal D}_1 F,$$
which follows from the definition of $\GG_2$, we find that
$$\|{\mathcal D}_1{\mathcal D}_0 \GG_2(F)\|_{L_1^2} \leq \|\partial_z^2 \GG_2(F)\|_{L_1^2} + \|{\mathcal D}_1 F\|_{L_1^2} \lesssim \|F_z\|_{L_1^2} + \|{\mathcal D}_1 F\|_{L_1^2}.\eqno{\Box}.$$

\begin{lemma} \label{Estimate for G3}
The estimate
$$\|\GG_3(\xi)\|_\star \lesssim \|\xi\|_{3/2}$$
holds for all $\xi \in H^{3/2}({\mathbb R})$.
\end{lemma}
{\bf Proof.} First note that
$$\tfrac{1}{2}|k|^2\left(\frac{I_0(|k|)^2}{I_1(|k|)^2}-1\right) \lesssim (1+|k|^2)^{1/2}$$
since
$$\frac{|k|^2}{2(1+|k|^2)^{1/2}}\left(\frac{I_0(|k|)^2}{I_1(|k|)^2}-1\right)
\rightarrow \left\{\begin{array}{ll}2 \quad \mbox{as $|k| \rightarrow 0$}, \\[2mm]
\frac{1}{2} \quad \mbox{as $|k| \rightarrow \infty$},\end{array}\right.$$
from which it follows that
\begin{align*}
\|\partial_z \GG_3(\xi)\|_{L_1^2}
& = \int_{-\infty}^\infty \int_0^1 r|k|^2\frac{I_0(|k|r)^2}{I_1(|k|)^2}|\hat{\xi}|^2 \dr\dk \\
& = \int_{-\infty}^\infty \int_0^1 \tfrac{1}{2}|k|^2\left(\frac{I_0(|k|)^2}{I_1(|k|)^2}-1\right)|\hat{\xi}|^2\dk\\
& \lesssim \|\xi\|_{1/2}^2
\end{align*}
and hence that
$$\|\partial_z^2 \GG_3(\xi)\|_{L_1^2} = \|\GG_3(\xi_z)\|_{L_1^2} \lesssim \|\xi_z\|_{1/2}^2.$$

Similarly
$$\tfrac{1}{2}|k|^2\left(1-\frac{I_0(|k|)^2}{I_1(|k|)^2}\right) +|k|\frac{I_0(|k|)}{I_1(|k|)}\lesssim (1+|k|^2)^{1/2}$$
since
$$\frac{|k|^2}{2(1+|k|^2)^{1/2}}\left(1-\frac{I_0(|k|)^2}{I_1(|k|)^2}\right)+\frac{|k|}{(1+|k|^2)^{1/2}}\frac{I_0(|k|)}{I_1(|k|)}
\rightarrow \left\{\begin{array}{ll}0 \quad \mbox{as $|k| \rightarrow 0$}, \\[2mm]
\frac{1}{2} \quad \mbox{as $|k| \rightarrow \infty$},\end{array}\right.$$
from which it follows that
\begin{align*}
\|{\mathcal D}_0 \GG_3(\xi)\|_{L_1^2}
& = \int_{-\infty}^\infty \int_0^1 r|k|^2 \frac{I_1(|k|r)^2}{I_1(|k|)^2}|\hat{\xi}|^2 \dr\dk \\
& = \int_{-\infty}^\infty \left(\tfrac{1}{2}|k|^2\left(1-\frac{I_0(|k|)^2}{I_1(|k|)^2}\right)+|k|\frac{I_0(|k|)}{I_1(|k|)}\right)|\hat{\xi}|^2 \dr\dk \\
& \lesssim \|\xi\|_{1/2}^2
\end{align*}
and hence that
$$\|\partial_z {\mathcal D}_0\GG_3(\xi)\|_{L_1^2} = \|{\mathcal D}_0\GG_3(\xi_z)\|_{L_1^2} \lesssim \|\xi_z\|_{1/2}^2.$$

Finally, using the identity
$${\mathcal D}_1{\mathcal D}_0 \GG_3(\xi) + \partial_z^2 \GG_3(\xi) = 0,$$
which follows from the definition of $\GG_3$, we find that
$$\|{\mathcal D}_1{\mathcal D}_0 \GG_3(\xi)\|_{L_1^2} \leq \|\partial_z^2 \GG_3(\xi)\|_{L_1^2} \lesssim \|\xi_z\|_{1/2}.\eqno{\Box}.$$

\subsection{Expansions}

Using the results in Section \ref{Everything is analytic} above, we obtain the expansions
\begin{equation}
u(\eta,\xi)=\sum_{j=0}^\infty u^j(\eta,\xi), \label{Expansion of u}
\end{equation}
where $u^j$ is homogeneous of degree $j$ in $\eta$ and linear in $\xi$, and
$$K(\eta)=\sum_{j=0}^\infty K_j(\eta), \qquad \KK(\eta)=\sum_{j=0}^\infty \KK_j(\eta), \qquad \LL(\eta)=\sum_{j=0}^\infty \LL_j(\eta),$$
where $K_j$, $\KK_j$, $\LL_j$ are homogeneous of degree $j$ in $\eta$. Note in particular the formulae
\begin{align}
\KK_1(\eta) & = (\gamma-1)\eta-\eta_{zz}, \label{KK1 formula v1} \\
\KK_2(\eta) & = A_0\eta^2 -\tfrac{1}{2}\eta_z^2, \label{KK2 formula v1} \\
\KK_3(\eta) & = B_0\eta^3 +\tfrac{1}{2}\eta\eta_z^2 +\tfrac{3}{2}\eta_z^2\eta_{zz}, \label{KK3 formula v1} \\
\nonumber \\
\LL_1(\eta) & = K_0 \eta, \label{LL1 formula v1} \\
\LL_2(\eta) & = \tfrac{1}{2}\big(\eta_z^2-(K_0\eta)^2+K_0\eta^2+2K_1(\eta)\eta\big), \label{LL2 formula v1} \\
\LL_3(\eta) & = -\eta_z^2K_0\eta-\tfrac{1}{2}(K_0\eta)(K_0\eta^2+2K_1(\eta)\eta)+\tfrac{1}{2}K_1(\eta)\eta^2+K_2(\eta)\eta, \label{LL3 formula v1}
\end{align}
where
$$A_0 = -\gamma-\tfrac{1}{2}\gamma\nu^{\prime\prime}(1)+1, \qquad B_0 = \gamma+\gamma\nu^{\prime\prime}(1)+\tfrac{1}{6}\gamma\nu^{\prime\prime\prime}(1)-1,$$
which are obtained from equations \eqref{Defn of KK}, \eqref{Defn of LL}.

The terms in the expansion of $u(\eta,\xi)$ can be computed by proceeding
formally. Substituting \eqref{Expansion of u}
into equations \eqref{Flattened BC 1}--\eqref{Flattened BC 2} and equating terms which are homogeneous of order $j$ in $\eta$ yields
a boundary-value problem for $u^j$ in terms of $u^0$, \ldots, $u^{j-1}$. Formulae for the terms in the expansion for $K(\eta)$ in terms of Fourier multipliers
are then recovered from the formula
$$K_j(\eta)(\xi) = -u_z^j(\eta,\xi)|_{r=1}.$$

\begin{remark} \label{Subtle cancellations}
This method leads to formulae involving ever more derivatives of $\eta$ and $\xi$ in the individual terms in the formulae for
$K_j(\eta)$; the overall validity of the formulae arises from subtle cancellations between the terms
(see Nicholls and Reitich \cite[\S2.2]{NichollsReitich01a} for a discussion of this phenomenon in the context of the classical Dirichlet--Neumann operator).
\end{remark}

\begin{proposition}
The operator $K_0 \in \LL(H^{3/2}({\mathbb R}),H^{1/2}({\mathbb R}))$ is given by the formula
$$K_0\xi = f(D)\xi,$$
which also defines an operator in ${\mathcal L}(H^{s+1}({\mathbb R}), H^s({\mathbb R}))$ for each $s \geq 0$.
\end{proposition}
{\bf Proof.} The solution to the boundary-value problem
\begin{alignat*}{2}
{\mathcal D}_1 {\mathcal D}_0 u^0 + u_{zz}^0 & = 0, \qquad && 0<r<1, \\
{\mathcal D}_0 u^0 & = \xi_z, && r=1,
\end{alignat*}
for $u^0$ is
$$u^0=\FF^{-1}\left[\frac{\i k I_0(|k|r)}{|k|I_1(|k|)}\hat{\xi}\right],$$
such that
$$K_0\xi = -(u_z^0)|_{r=1}=f(D)\xi.$$
Furthermore, the estimate
$$f(k) \lesssim \sqrt{1+k^2}, \qquad k \in {\mathbb R},$$
which follows from the calculation
$$\frac{f(k)}{\sqrt{1+k^2}} = \frac{|k|I_0(|k|)}{\sqrt{1+k^2}I_1(|k|)}
\rightarrow \left\{\begin{array}{ll}2 \quad \mbox{as $|k| \rightarrow 0$}, \\[2mm]
1 \quad \mbox{as $|k| \rightarrow \infty$},\end{array}\right.$$
implies that
$$\|K_0 \xi\|_s \lesssim \|\xi\|_{s+1}$$
for all $s \geq 0$.\qed

\begin{proposition}$ $
The operators $K_1$ and $K_2$ are given by the formulae
\begin{align*}
K_1(\eta)\xi &= -(\eta\xi_z)_z - K_0(\eta K_0\xi), \\
K_2(\eta)\xi & = \tfrac{1}{2}(\eta^2 K_0\xi)_{zz} + \tfrac{1}{2}K_0(\eta^2 \xi_{zz}) + \tfrac{1}{2}(\eta^2\xi_z)_z
-\tfrac{1}{2}K_0(\eta^2K_0\xi) + K_0(\eta K_0(\eta K_0\xi))
\end{align*}
for each $\eta \in H^2({\mathbb R})$ and $\xi \in H^{3/2}({\mathbb R})$ (see Remark \ref{Subtle cancellations}).
\end{proposition}
{\bf Proof.} The solution to the boundary-value problem
\begin{alignat*}{2}
{\mathcal D}_1 {\mathcal D}_0 u^1 + u_{zz}^1 & = {\mathcal D_1}(r\eta_zu_z^0) + \partial_z(r\eta_z{\mathcal D}_0 u^0-2\eta u_z^0), \qquad\qquad && 0<r<1, \\
{\mathcal D}_0 u^1& = r\eta_zu_z^0, && r=1,
\end{alignat*}
for $u^1$ is
$$u^1=r\eta {\mathcal D}_0 u^0 + w^1,$$
where
\begin{alignat*}{2}
{\mathcal D}_1 {\mathcal D}_0 w^1 + w_{zz}^1 & = 0, \qquad\qquad && 0<r<1, \\
{\mathcal D}_0 w^1 & = (\eta u_z^0)_z, && r=1,
\end{alignat*}
such that
$$u^1=r\eta {\mathcal D}_0 u^0 + \FF^{-1}\left[\frac{\i kI_0(|k|r)}{ |k|I_1(|k|)}\FF[\eta u_z^0|_{r=1}]\right].$$
It follows that
$$K_1(\eta)\xi = -(u_z^1)|_{r=1}=-(\eta\xi_z)_z-K_0(\eta K_0\xi).$$

Similarly, the solution to the boundary-value problem
\begin{alignat*}{2}
{\mathcal D}_1 {\mathcal D}_0 u^2 + u_{zz}^2 & = {\mathcal D}_1(r\eta_zu_z^1+r\eta\eta_zu_z^0-r^2\eta_z{\mathcal D}_0 u^0), \qquad\qquad\qquad && 0<r<1, \\
& \qquad\mbox{}+\partial_z (r\eta_z {\mathcal D}_0 u^1 + r\eta\eta_z {\mathcal D}_0 u^0-\eta^2 u_z^0-2\eta u_z^1), \\
{\mathcal D}_0 u^2 & = r\eta_zu_z^1+r\eta\eta_zu_z^0-r^2\eta_z{\mathcal D}_0 u^0, && r=1,
\end{alignat*}
for $u^2$ is
$$u^2=-\tfrac{1}{2}\eta^2{\mathcal D}_0 (r^2{\mathcal D}_0 u^0)+r\eta {\mathcal D}_0 u^1+w^2,$$
where
\begin{alignat*}{2}
{\mathcal D}_1 {\mathcal D}_0 w^2 + w_{zz}^2 & = 0, && 0<r<1, \\
{\mathcal D}_0 w^2 & = (u_z^1)_z + (\tfrac{1}{2}\eta^2 u_z^0)_z - (\tfrac{1}{2}(\eta^2 {\mathcal D}_0 u^0)_z)_z, \qquad && r=1,
\end{alignat*}
such that
\begin{align*}
u^2 & = -\tfrac{1}{2}\eta^2{\mathcal D}_0 (r^2{\mathcal D}_0 u^0) + r\eta {\mathcal D}_0 u^1 
+\FF^{-1}\left[\frac{\i kI_0(|k|r)}{ |k|I_1(|k|)} \FF[\tfrac{1}{2}\eta^2 u_z^0|_{r=1}]\right] \\
& \qquad\mbox{}-\FF^{-1}\left[\frac{\i  kI_0(|k|r)}{ |k|I_1(|k|)} \FF[\tfrac{1}{2}(\eta^2 {\mathcal D}_0 u^0)_z|_{r=1}]\right]
+ \FF^{-1}\left[\frac{\i  kI_0(|k|r)}{ |k|I_1(|k|)} \FF[\eta u_z^1|_{r=1}]\right].
\end{align*}
We conclude that
\begin{align*}
K_2(\eta)\xi & = -(u_z^2)|_{r=1} \\
& = \tfrac{1}{2}(\eta^2 \mathcal{D}_0^2 u^0)_z + (\eta^2 {\mathcal D}_0 u^0)_z - (\eta \eta_z u_z^0)_z - \tfrac{1}{2}K_0((\eta^2\xi_z)_z)
+\tfrac{1}{2}K_0(\eta^2 u_z^0) - K_0(\eta K_1(\eta)\xi) \\
& = \tfrac{1}{2}(\eta^2 K_0\xi)_{zz} + \tfrac{1}{2}K_0(\eta^2 \xi_{zz}) + \tfrac{1}{2}(\eta^2\xi_z)_z
-\tfrac{1}{2}K_0(\eta^2K_0\xi) + K_0(\eta K_0(\eta K_0\xi)).\alignqed
\end{align*}
\begin{corollary} \label{Formulae for 123}
The formulae
\begin{align*}
\KK_1(\eta) & = (\gamma-1)\eta-\eta_{zz}, \\
\KK_2(\eta) & = A_0\eta^2 -\tfrac{1}{2}\eta_z^2, \\
\KK_3(\eta) & = B_0\eta^3 +\tfrac{1}{2}\eta\eta_z^2 +\tfrac{3}{2}\eta_z^2\eta_{zz}, \\
\\
\LL_1(\eta) & = K_0 \eta, \\
\LL_2(\eta) & = \tfrac{1}{2}(\eta_z^2-(K_0\eta)^2-(\eta^2)_{zz}-2K_0(\eta K_0\eta) + K_0\eta^2), \\
\LL_3(\eta) & = \tfrac{1}{2}(K_0\eta)(\eta^2)_{zz} + (K_0\eta)(K_0(\eta K_0\eta)) - \tfrac{1}{2}(K_0\eta)(K_0\eta^2)-\eta_z^2(K_0\eta) \\
& \qquad \mbox{}+\tfrac{1}{2}(\eta^2K_0 \eta)_{zz}+\tfrac{1}{2}K_0(\eta^2\eta_{zz})-\tfrac{1}{2}(\eta^2\eta_z)_z -\tfrac{1}{2}K_0(\eta^2K_0\eta) \\
& \qquad \mbox{}+K_0(\eta K_0(\eta K_0\eta))-\tfrac{1}{2}K_0(\eta K_0 \eta^2)
\end{align*}
hold for all $\eta \in H^2({\mathbb R}^2)$ (see Remark \ref{Subtle cancellations}).
\end{corollary}

%

\begin{lemma} $ $ \label{2, 3 and r estimates}
\begin{itemize}
\item[(i)]
The estimates
\begin{align*}
\|\KK_2(\eta)\|_0 \lesssim \|\eta\|_\ZZ \|\eta\|_2, \qquad & \|\mathrm{d}\KK_2[\eta](\rho)\|_0 \lesssim \|\eta\|_\ZZ \|\rho\|_2, \\
\|\LL_2(\eta)\|_0 \lesssim \|\eta\|_\ZZ \|\eta\|_2, \qquad & \|\mathrm{d}\LL_2[\eta](\rho)\|_0 \lesssim \|\eta\|_\ZZ \|\rho\|_2
\end{align*}
hold for all $\eta$, $\rho \in H^2({\mathbb R})$.
\item[(ii)]
The estimates
\begin{align*}
\|\KK_3(\eta)\|_0 \lesssim \|\eta\|_\ZZ^2 \|\eta\|_2, \qquad & \|\mathrm{d}\KK_3[\eta](\rho)\|_0 \lesssim \|\eta\|_\ZZ^2 \|\rho\|_2 + \|\eta\|_\ZZ \|\eta\|_2 \|\rho\|_2, \\
\|\LL_3(\eta)\|_0 \lesssim \|\eta\|_\ZZ^2 \|\eta\|_2, \qquad & \|\mathrm{d}\LL_3[\eta](\rho)\|_0 \lesssim \|\eta\|_\ZZ^2 \|\rho\|_2 + \|\eta\|_\ZZ \|\eta\|_2 \|\rho\|_2
\end{align*}
hold for all $\eta \in U$ and $\rho \in H^2({\mathbb R})$.
\item[(iii)]
The quantities
$$\KK_\mathrm{r}(\eta) = \sum_{j=4}^\infty \KK_j(\eta), \qquad \LL_\mathrm{r}(\eta) = \sum_{j=4}^\infty \LL_j(\eta)$$
satisfy the estimates
\begin{align*}
\|\KK_\mathrm{r}(\eta)\|_0 \lesssim \|\eta\|_\ZZ^3 \|\eta\|_2, \qquad & \|\mathrm{d}\KK_r[\eta](\rho)\|_0 \lesssim \|\eta\|_\ZZ^3 \|\rho\|_2 + \|\eta\|_\ZZ^2 \|\eta\|_2 \|\rho\|_2, \\
\|\LL_\mathrm{r}(\eta)\|_0 \lesssim \|\eta\|_\ZZ^3 \|\eta\|_2, \qquad & \|\mathrm{d}\LL_r[\eta](\rho)\|_0 \lesssim \|\eta\|_\ZZ^3 \|\rho\|_2 + \|\eta\|_\ZZ^2 \|\eta\|_2 \|\rho\|_2
\end{align*}
hold for all $\eta \in U$ and $\rho \in H^2({\mathbb R})$.
\end{itemize}
\end{lemma}
{\bf Proof.} These results are obtained by estimating the right-hand sides of \eqref{KK2 formula v1}, \eqref{KK3 formula v1}, \eqref{LL2 formula v1}, \eqref{LL3 formula v1} and
\begin{align*}
\KK_\mathrm{r}(\eta) & = \KK(\eta) - \KK_1(\eta) - \KK_2(\eta) - \KK_3(\eta) \\
& = - \left(\gamma \nu\left(\frac{1}{1+\eta}\right)-\gamma \nu(1) + \gamma \eta + (A_0-1)\eta^2 + (B_0+1)\eta^3\right) \\
& \qquad\mbox{}-\left(\frac{1}{(1+\eta_z^2)^{3/2}}-1+\tfrac{3}{2}\eta_z^2\right)\eta_{zz} + (1-\eta)\left(\frac{1}{(1+\eta_z^2)^{1/2}}-1+\tfrac{1}{2}\eta_z^2\right) \\
& \qquad\mbox{}+(\eta^2-\eta^3)\left(\frac{1}{(1+\eta_z^2)^{1/2}}-1\right)+\left(\frac{1}{1+\eta}-1+\eta-\eta^2+\eta^3\right)\frac{1}{(1+\eta_z^2)^{1/2}}, \\
\\
\LL_\mathrm{r}(\eta) & = \LL(\eta) - \LL_1(\eta) - \LL_2(\eta) - \LL_3(\eta) \\
& = -\tfrac{1}{2}(K_{\geq 2}(\eta)\eta + \tfrac{1}{2}K(\eta)\eta^2)^2 - (K_0\eta+K_1(\eta)\eta)(K_{\geq 2}(\eta)\eta+\tfrac{1}{2}K_{\geq 1}(\eta)\eta^2)-\tfrac{1}{2}\left(K_1(\eta)\eta\right)^2\\
& \qquad\mbox{}-\tfrac{1}{2}K_1(\eta)\eta(K_0\eta^2)+\frac{1}{2}\left(\frac{1}{1+\eta_z^2}-1\right)(\eta_z - \eta_zK(\eta)\eta-\tfrac{1}{2}\eta_zK(\eta)\eta^2)^2 \\
& \qquad\mbox{}+\tfrac{1}{2}(\eta_zK(\eta)\eta+\tfrac{1}{2}\eta_zK(\eta)\eta^2)^2-\eta_z^2(K_{\geq 1}(\eta)\eta+\tfrac{1}{2}K(\eta)\eta^2) \\
& \qquad\mbox{}+K_{\geq 3}(\eta)\eta+\tfrac{1}{2}K_{\geq 2}(\eta)\eta^2
\end{align*}
with
$$K_{\geq 1}(\eta) = \sum_{j=1}^\infty K_j(\eta), \qquad K_{\geq 2}(\eta) = \sum_{j=2}^\infty K_j(\eta), \qquad K_{\geq 3}(\eta) = \sum_{j=3}^\infty K_j(\eta)$$
using the methods described in the proof of Corollary \ref{KK and LL are analytic}, noting that
$$\|K_j(\eta)\eta\|_{1/2} \lesssim \|\eta\|_\ZZ^j \|\eta\|_{3/2}, \qquad \|K(\eta)\eta\|_{1/2} \lesssim \|\eta\|_{3/2}$$
and
$$\|K_{\geq 1}(\eta)\eta\|_{1/2} \lesssim \|\eta\|_\ZZ\|\eta\|_{3/2}, \qquad
\|K_{\geq 2}(\eta)\eta\|_{1/2} \lesssim \|\eta\|_\ZZ^2 \|\eta\|_{3/2}, \qquad
\|K_{\geq 3}(\eta)\eta\|_{1/2} \lesssim \|\eta\|_\ZZ^3 \|\eta\|_{3/2}.$$

The estimates for the derivatives are obtained in the same way.\qed

\section{Reduction} \label{Reduction}

In this section we reduce the equation
\begin{equation}
{\mathcal K}(\eta)-c_0^2(1-\varepsilon^2){\mathcal L}(\eta)=0 \label{Basic equation}
\end{equation}
to a perturbation of a full-dispersion model equation using a technique reminiscent of the Lyapunov-Schmidt reduction. We work
in the subset $U$ of the basic space $\XX=H^2({\mathbb R}^2)$ (see equation \eqref{Defn of U}, so that equation \eqref{Basic equation} holds in $L^2({\mathbb R}^2)$.
Respecting the decomposition of $\eta$ into two parts, we decompose
$\XX$ into the direct sum of the spaces
$$\XX_1 = \chi(D)\XX, \qquad \XX_2 = (1-\chi(D))\XX$$
and equip $\XX_1$ and $\XX_2$ with respectively the scaled norm
$$
\nn \eta_1 \nn^2 = \int_{\mathbb R} (1+\varepsilon^{-2}(|k|-\omega)^2)|\hat{\eta}_1|^2\dk$$
(with the convention that $\omega=0$ if $1<\gamma<9$)
and the usual norm for $H^2({\mathbb R})$.

\begin{proposition} \label{Scaled space estimates}
The estimate
$$\|\hat{\eta}\|_{L^1{\mathbb R})} \lesssim \varepsilon^{1/2}\nn\eta_1\nn$$
holds for every $\eta \in \XX_1$, and in particular
$$\|\eta\|_\ZZ \lesssim \varepsilon^{1/2}\nn \eta_1 \nn + \|\eta_2\|_2$$
for every $\eta \in \XX$.
\end{proposition}
{\bf Proof.} This result follows from the calculation
\begin{align*}
\int_{\mathbb R} |\hat{\eta}_1(k)|\dk &= \int_{\mathbb R} \frac{(1\!+\!\varepsilon^{-2}(|k|\!-\!\omega)^2)^{1/2}}{(1\!+\!\varepsilon^{-2}(|k|\!-\!\omega)^2)^{1/2}}|\hat{\eta}_1(k)|\dk\\
&
\leq \left(\int_{\mathbb R} \frac{1}{1\!+\!\varepsilon^{-2}(|k|\!-\!\omega)^2}\dk\right)^{\!\!1/2}\!\!\!\!\nn \eta_1 \nn\\
& = \left(\pi \varepsilon+2\varepsilon\arctan\frac{\omega}{\varepsilon}\right)^{1/2} \nn \eta_1 \nn.\alignqed
\end{align*}

Clearly $\eta \in U$ satisfies \eqref{Basic equation} if and only if
\begin{align*}
\chi(D)\left(\KK(\eta_1+\eta_2)-c_0^2(1-\varepsilon^2)\LL(\eta_1+\eta_2)\right)&=0, \\
(1-\chi(D))\left(\KK(\eta_1+\eta_2)-c_0^2(1-\varepsilon^2)\LL(\eta_1+\eta_2)\right)&=0,
\end{align*}
and these equations can be rewritten as
\begin{align}
g(D)\eta_1
+c_0^2\varepsilon^2K_0\eta_1+\chi(D)\NN(\eta_1+\eta_2)&=0, \label{X1 eqn} \\
g(D)\eta_2 +c_0^2\varepsilon^2 K_0\eta_2+(1-\chi(D))\NN(\eta_1+\eta_2)&=0,
\label{X2 eqn}
\end{align}
in which
$$\NN(\eta)=\KK_2(\eta)+\KK_3(\eta)+\KK_\mathrm{r}(\eta)-c_0^2(1-\varepsilon^2)(\LL_2(\eta)+\LL_3(\eta)+\LL_\mathrm{r}(\eta)).$$

We proceed by writing \eqref{X2 eqn} as a fixed-point equation for $\eta_2$ using Proposition \ref{invert g},
which follows from the fact that $g(k) \gtrsim |k|^2$ for $k \not\in S$,
and solving it for $\eta_2$ as a function of $\eta_1$ using Theorem \ref{thm:fixed-point},
which is proved by a straightforward application of the contraction mapping principle. Substituting $\eta_2=\eta_2(\eta_1)$ into
\eqref{X1 eqn} yields a reduced equation for $\eta_1$, which can be rewritten as a perturbation of a full-dispersion model equation by
applying a further change of variable. Full details are given in Sections \ref{Reduction - sst} and \ref{Reduction - wst} below, which
deal with the cases $1 < \gamma < 9$ (`strong surface tension') and $\gamma>9$ (`weak surface tension') separately.

\begin{proposition} \label{invert g}
The mapping $(1-\chi(D))g(D)^{-1}$ is a bounded linear operator $L^2({\mathbb R}) \rightarrow \XX_2$.
\end{proposition}
\begin{theorem}
\label{thm:fixed-point}
Let $\XX_1$, $\XX_2$ be Banach spaces, $X_1$, $X_2$ be closed, convex sets in, respectively, $\XX_1$, $\XX_2$ containing the origin and $\GG\colon X_1\times X_2 \to \XX_2$ be a smooth function. Suppose that there exists a 
continuous function $r\colon X_1\to [0,\infty)$ such that
$$
\|\GG(x_1,0)\|\le \tfrac{1}{2}r, \quad \|\mathrm{d}_2 \GG[x_1,x_2]\|\le \tfrac{1}{3}
$$ 
for each $x_2\in \bar B_r(0)\subseteq X_2$ and each $x_1\in X_1$.

Under these hypotheses there exists for each $x_1\in X_1$ a unique solution $x_2=x_2(x_1)$ of the fixed-point equation
$x_2=\GG(x_1,x_2)$
satisfying $x_2(x_1)\in \bar B_r(0)$. Moreover $x_2(x_1)$ is a smooth function of $x_1\in X_1$ and in particular satisfies the estimate
$$
\|\mathrm{d} x_2[x_1]\|\le 2\|\mathrm{d}_1 \GG[x_1, x_2(x_1)]\|.
$$
\end{theorem}

\subsection{Strong surface tension} \label{Reduction - sst}

Suppose that $1 < \gamma < 9$. We write \eqref{X2 eqn} in the form
\begin{equation}
\eta_2 = -(1-\chi(D))g(D)^{-1}\left(c_0^2\varepsilon^2 K_0\eta_2+\NN(\eta_1+\eta_2)\right)
\label{sst eta2 eqn}
\end{equation}
and apply Theorem \ref{thm:fixed-point} with
$$
X_1=\{\eta_1\in \XX_1 \colon \nn \eta_1\nn \le R_1\}, \qquad
X_2=\{\eta_2\in \XX_2 \colon \| \eta_2\|_2 \le R_2\};
$$
the function $\GG$ is given by the right-hand side of \eqref{sst eta2 eqn}. 
Using Proposition \ref{Scaled space estimates}
one can guarantee that\linebreak $\|\hat{\eta}_1\|_{L^1({\mathbb R}^2)} < \frac{1}{2}M$ for all
$\eta_1 \in X_1$ for an arbitrarily large value
of $R_1$; the value of $R_2$ is
constrained by the requirement that $\|\eta_2\|_2 < \frac{1}{2}M$ for all $\eta_2 \in X_2$. The next lemma follows from Lemma \ref{2, 3 and r estimates}, its corollary from Proposition \ref{invert g}.

\begin{lemma}
The estimates
\begin{itemize}
\item[(i)]
$\|\NN(\eta_1,\eta_2)\|_0\lesssim \varepsilon^{1/2} \nn \eta_1\nn^2+\varepsilon^{1/2} \nn \eta_1\nn \|\eta_2\|_2
+\nn \eta_1 \nn \|\eta_2\|_2^2 +\|\eta_2\|_2^2$,
\item[(ii)]
$\|\mathrm{d}_1\NN[\eta_1,\eta_2]\|_{\LL(\XX_1,L^2({\mathbb R}))}\lesssim \varepsilon^{1/2} \nn \eta_1\nn+\varepsilon^{1/2}\|\eta_2\|_2+ \|\eta_2\|_2^2$,
\item[(iii)]
$\|\mathrm{d}_2\NN[\eta_1,\eta_2]\|_{\LL(\XX_2,L^2({\mathbb R}))}\lesssim \varepsilon^{1/2} \nn \eta_1 \nn + \nn \eta_1 \nn \|\eta_2\|_2 + \|\eta_2\|_2$,
\end{itemize}
where with a slight abuse of notation we write $\NN(\eta_1+\eta_2)$ as $\NN(\eta_1,\eta_2)$,
hold for each $\eta_1\in X_1$ and $\eta_2\in X_2$.
\end{lemma}
\begin{corollary}
\label{Complete sst GG estimates}
The estimates
\begin{itemize}
\item[(i)]
$\|\GG(\eta_1,\eta_2)\|_2\lesssim \varepsilon^{1/2} \nn \eta_1\nn^2+\varepsilon^{1/2} \nn \eta_1\nn \|\eta_2\|_2
+\nn \eta_1 \nn \|\eta_2\|_2^2 +\|\eta_2\|_2^2 + \varepsilon^2 \|\eta_2\|_2$,
\item[(ii)]
$\|\mathrm{d}_1\GG[\eta_1,\eta_2]\|_{\LL(\XX_1,\XX_2)}\lesssim \varepsilon^{1/2} \nn \eta_1\nn+\varepsilon^{1/2}\|\eta_2\|_2+ \|\eta_2\|_2^2$,
\item[(iii)]
$\|\mathrm{d}_2\GG[\eta_1,\eta_2]\|_{\LL(\XX_1,\XX_2)}\lesssim \varepsilon^{1/2} \nn \eta_1 \nn + \nn \eta_1 \nn \|\eta_2\|_2 + \|\eta_2\|_2+\varepsilon^2$
\end{itemize}
hold for each $\eta_1\in X_1$ and $\eta_2\in X_2$.
\end{corollary}

\begin{theorem} \label{thm:estimate eta2}
Equation \eqref{sst eta2 eqn} has a unique solution $\eta_2 \in
X_2$ which depends smoothly upon $\eta_1 \in X_1$ and satisfies the estimates
$$
\|\eta_2(\eta_1)\|_2 \lesssim \varepsilon^{1/2} \nn \eta_1\nn^2,  \qquad
\|\mathrm{d}\eta_2[\eta_1]\|_{\LL(\XX_1,\XX_2)} \lesssim \varepsilon^{1/2} \nn \eta_1\nn.
$$
\end{theorem}
{\bf Proof.}
Choosing $R_2$ and $\varepsilon$ sufficiently small and setting
$r(\eta_1)=\sigma \varepsilon^{1/2} \nn \eta_1 \nn^2$ for a sufficiently large value
of $\sigma>0$, one finds that
\[\|\GG(\eta_1,0)\|_2 \lesssim \tfrac{1}{2}r(\eta_1), \qquad
\|\mathrm{d}_2 \GG[\eta_1,\eta_3]\|_{\LL(\XX_2,\XX_2)}  \lesssim \varepsilon^{1/2}
\]
for $\eta_1 \in X_1$ and $\eta_2 \in \overline{B}_{r(\eta_1)}(0) \subset X_2$
(Corollary \ref{Complete sst GG estimates}(i), (iii)).
Theorem \ref{thm:fixed-point} asserts that equation \eqref{sst eta2 eqn} has a unique solution $\eta_2$
in $\overline{B}_{r(\eta_1)}(0) \subset X_2$ which depends smoothly upon $\eta_1 \in X_1$, and
the estimate for its derivative follows from Corollary \ref{Complete sst GG estimates}(ii).\qed\medskip

Substituting $\eta_2=\eta_2(\eta_1)$ into \eqref{X1 eqn} yields the reduced equation
\begin{align}
g(D)\eta_1 +c_0^2\varepsilon^2K_0\eta_1+\chi(D)\NN(\eta_1+\eta_2(\eta_1))=0
\label{sst red eq v1}
\end{align}
for $\eta_1 \in X_1$. The leading-order terms in this equation are computed by approximating the operators
$\partial_z$ and $K_0$ in its quadratic part by constants.

\begin{proposition} \label{sst approximate identities}
The estimates
\begin{itemize}
\item[(i)]
$\eta_{1z} = O(\varepsilon \nn \eta_1\nn)$,
\item[(ii)]
$K_0 \eta_1 = 2\eta_1 + O(\varepsilon\nn \eta_1\nn)$,
\item[(iii)]
$K_0( \eta_1\rho_1) = 2\eta_1\rho_1 + O(\varepsilon^{3/2}\nn \eta_1\nn\nn \rho_1\nn)$
\end{itemize}
hold for all $\eta_1$, $\rho_1 \in \XX_1$. The order-of-magnitude estimates are computed with respect to the $L^2({\mathbb R})$-norm (which is equivalent to the
$H^s({\mathbb R})$-norm on the space $\chi(D)H^s({\mathbb R})$ for any $s \geq 0$).
\end{proposition}
{\bf Proof.}
This result follows from the calculations
\begin{align*}
\|\eta_{1z}\|_0 &= \||k| \hat{\eta}_1\|_0 \leq \varepsilon \nn \eta_1 \nn, \\[2mm]
\|(K_0-2I)\eta_1\|_0
& = \left\|(f(k)-2)\hat{\eta}_1\right\|_0
\lesssim \| |k| \hat{\eta}_1\|_0
\leq \varepsilon \nn \eta_1 \nn, \\[2mm]
\|(K_0-2I)(\eta_1\rho_1)\|_0
&\lesssim \left\||k|
\int_{{\mathbb R}} | \hat \eta_1(k-\tilde{k})| | \hat \rho_1(\tilde{k})| \dtildek\right\|_0\\
&\lesssim \left\|\int_{{\mathbb R}}  |k-\tilde{k}|  | \hat \eta_1(k-\tilde{k})|  | \hat \rho_1(\tilde{k})| \dtildek+\int_{{\mathbb R}} |\tilde{k}| | \hat \eta_1(k-\tilde{k})| | \hat \rho_1(\tilde{k})| \dtildek\right\|_0 \\
&\lesssim \|\ |k|\hat \eta_1 \|_0 \|\hat \rho_1\|_{L^1({\mathbb R})}+\|\hat \eta_1\|_{L^1({\mathbb R})}
\|\ |k|\hat \rho_1\|_0 \\
&\lesssim \varepsilon^{3/2}\nn \eta_1 \nn\,  \nn \rho_1\nn
\end{align*}
for each $\eta_1$, $\rho_1 \in X_1$, where we have also used Young's inequality.\qed

The leading-order terms in the nonlinear part of \eqref{sst red eq v1} are now obtained from Corollary \ref{First rf estimate} (which follows from Corollary \ref{Formulae for 123} and
Proposition \ref{sst approximate identities}) and Lemma \ref{Second rf estimate} (which follows from Lemma \ref{2, 3 and r estimates}) below.
Here we use the symbol $\underline{O}(\varepsilon^s \nn \eta_1\nn^t)$ (with $s \geq 0$, $t \geq 1$)
to denote a smooth function
 $\RR^\varepsilon: X_1 \rightarrow L^2({\mathbb R})$ which satisfies the estimates 
$$
\|\RR^\varepsilon(\eta_1)\|_0 \lesssim \varepsilon^s  \nn \eta_1\nn^t \quad
\|\mathrm{d}\RR^\varepsilon[\eta_1]\|_{\LL(\XX_1,L^2({\mathbb R}^2))}\lesssim \varepsilon^s  \nn \eta_1 \nn^{t-1}$$
for each $\eta_1 \in X_1$.
\begin{corollary} \label{First rf estimate}
The estimates
\begin{itemize}
\item[(i)]
$\KK_2(\eta_1+\eta_2(\eta_1))=\left(-\gamma -\tfrac{1}{2}\gamma v^{\prime\prime}(1)+1\right)\eta_1^2 + \underline{O}(\varepsilon \nn \eta_1\nn^2)$,
\item[(ii)]
$\LL_2(\eta_1+\eta_2(\eta_1))=-5\eta_1^2 + \underline{O}(\varepsilon \nn \eta_1\nn^2)$
\end{itemize}
hold for each $\eta_1 \in X_1$.
\end{corollary}
\begin{lemma} \label{Second rf estimate}
The estimate
$$
\NN(\eta_1+\eta_2(\eta_1)) = \KK_2(\eta_1)-c_0^2(1-\varepsilon^2)\LL_2(\eta_1)+\underline{O}(\varepsilon\nn \eta_1 \nn^3)
$$
holds for each $\eta_1 \in X_1$.
\end{lemma}

We conclude that the reduced equation for $\eta_1$ is the \emph{perturbed full dispersion Korteweg-de Vries equation}
$$g(D)\eta_1 + c_0^2\varepsilon^2 K_0\eta_1
+\chi(D)\Big(2c_0^2d_0\eta_1^2 + \underline{O}(\varepsilon \nn \eta_1 \nn^2)\Big)=0,$$
and applying Proposition \ref{sst approximate identities}(ii), one can further simplify it to
$$g(D)\eta_1 + 2c_0^2\varepsilon^2 \eta_1 +\chi(D)\Big(2c_0^2d_0\eta_1^2 + \underline{O}(\varepsilon \nn \eta_1 \nn^2)
+ \underline{O}(\varepsilon^3 \nn \eta_1 \nn)\Big)=0.$$
Finally, we introduce the Korteweg-de Vries scaling
$$\eta_1(z)=\varepsilon^2 \zeta(\varepsilon z),$$
noting that $I: \eta_1 \rightarrow \zeta$ is an isomorphism $\XX_1 \rightarrow H_\varepsilon^1({\mathbb R})$ and
$\chi(D)L^2({\mathbb R}) \rightarrow L_\varepsilon^2({\mathbb R})$ and choosing $R>1$ large enough so that $\zeta_\mathrm{KdV}\in B_R(0)$
(and $\varepsilon>0$ small enough so that $B_R(0) \subseteq H_\varepsilon^1({\mathbb R})$ is contained in $I[X_1]$). We find that
$\zeta \in B_R(0) \subseteq H_\varepsilon^1({\mathbb R})$ satisfies the equation
\begin{equation}
\varepsilon^{-2}g(\varepsilon D)\zeta + 2c_0^2\zeta +2c_0^2d_0\chi_0(\varepsilon D) \zeta^2 + \varepsilon^{1/2}\underline{O}^\varepsilon_0(\| \zeta \|_1)=0,
\label{PFDKdV - freshly derived}
\end{equation}
which holds in $L_\varepsilon^2({\mathbb R})$, where the symbol $D$ now means $-\i \partial_Z$ and
the symbol $\underline{O}^\varepsilon_n(\varepsilon^s\| \zeta \|_1^t)$
denotes a smooth function
$\RR: B_R(0) \subseteq H_\varepsilon^1({\mathbb R}) \rightarrow H_\varepsilon^n({\mathbb R})$
which satisfies the estimates
$$
\|\RR(\zeta)\|_n \lesssim \varepsilon^s\| \zeta \|_1^t \quad
\|\mathrm{d}\RR[\zeta]\|_{\LL(H^1({\mathbb R}),H^n({\mathbb R}))}\lesssim \varepsilon^s\| \zeta \|_1^{t-1}
$$
for each $\zeta \in B_R(0)\subseteq H_\varepsilon^1({\mathbb R})$  (with $t\geq 1$, $s$, $n \geq 0$).
Note that $\nn \eta \nn = \varepsilon^{3/2} \|\zeta\|_1$ and the change of variable from $z$ to $Z=\varepsilon z$ introduces
an additional factor of $\varepsilon^{1/2}$ in the remainder term.

Equation \eqref{sst red eq v1} is invariant under the reflection $\eta_1(z) \mapsto \eta_1(-z)$;
a familiar argument shows that it is inherited from the corresponding invariance of \eqref{X1 eqn}, \eqref{sst eta2 eqn}
under $\eta_1(z) \mapsto \eta_1(-z)$,
$\eta_2(z) \mapsto \eta_2(-z)$  when applying Theorem \ref{thm:fixed-point}.
The invariance is likewise inherited by \eqref{PFDKdV - freshly derived}, which is invariant under the reflection $\zeta(Z) \mapsto \zeta(-Z)$.

\subsection{Weak surface tension}  \label{Reduction - wst}

Suppose that $\gamma>9$. Since $\chi(D)\KK_2(\eta_1)$ and $\chi(D)\LL_2(\eta_1)$ both vanish the nonlinear term in
\eqref{X1 eqn} is at leading order cubic in $\eta_1$, so that this equation may be rewritten as
\begin{equation}
g(D)\eta_1+c_0^2\varepsilon^2K_0\eta_1+\chi(D)\left(\NN(\eta_1+\eta_2)+c_0^2(1-\varepsilon^2)\LL_2(\eta_1)-\KK_2(\eta_1)\right)=0.
\label{cubic eta1 eqn}
\end{equation}
To compute the reduced equation for $\eta_1$ we need an explicit formula for the leading-order quadratic part of $\eta_2(\eta_1)$; inspecting
\eqref{X2 eqn} shows that it is given by
\begin{equation}
F(\eta_1):=
(1-\chi(D))g(D)^{-1}\left(c_0^2(1-\varepsilon^2)\LL_2(\eta_1)-\KK_2(\eta_1)\right),
\label{Defn of F}
\end{equation}
an estimate for which is found using Lemma \ref{2, 3 and r estimates} (note that $K_0F(\eta)$ satisfies the same estimates as $F(\eta)$ since $\FF[F(\eta)]$ has
compact support).

\begin{proposition} \label{prop:F estimates}
The estimates
\begin{itemize}
\item[(i)]
$\|F(\eta_1)\|_2, \|K_0F(\eta_1)\|_2\lesssim \varepsilon^{1/2} \nn \eta_1\nn^2$,
\item[(ii)]
$\|\mathrm{d}F[\eta_1]\|_{\LL(\XX_1, \XX_2)},\|\mathrm{d}K_0F[\eta_1]\|_{\LL(\XX_1, \XX_2)}   \lesssim \varepsilon^{1/2} \nn \eta_1\nn$
\end{itemize}
hold for each $\eta_1\in X_1$.
\end{proposition}

It is convenient to write $\eta_2=F(\eta_1)+\eta_3$ and \eqref{X2 eqn}
in the form
\begin{equation}
\eta_3 = -(1-\chi(D))g(D)^{-1}\Big(\NN(\eta_1+F(\eta_1)+\eta_3)+c_0^2(1-\varepsilon^2)\LL_2(\eta_1)-\KK_2(\eta_1)
+c_0^2\varepsilon^2 K_0(F(\eta_1)+\eta_3)\Big) \label{wst eta3 eqn}
\end{equation}
(with the requirement that $\eta_1+F(\eta_1)+\eta_3 \in U$).
We apply Theorem \ref{thm:fixed-point} to equation \eqref{wst eta3 eqn} with
$$
X_1=\{\eta_1\in \XX_1 \colon \nn \eta_1\nn \le R_1\}, \qquad
X_3=\{\eta_3\in \XX_2 \colon \| \eta_3\|_3 \le R_3\};
$$
the function $\GG$ is given by the right-hand side of \eqref{wst eta3 eqn}. (Here we write
$X_3$ rather than $X_2$ for notational clarity.) Using Proposition \ref{Scaled space estimates}
one can guarantee that $\|\hat{\eta}_1\|_{L^1({\mathbb R})} < \frac{1}{2}M$ for all
$\eta_1 \in X_1$ for an arbitrarily large value
of $R_1$; the value of $R_3$
is constrained by the requirement that $\|F(\eta_1) + \eta_3\|_2 < \frac{1}{2}M$ for all $\eta_1 \in X_1$ and $\eta_3 \in X_3$,
so that $\eta_1+F(\eta_1)+\eta_3 \in U$
(Proposition \ref{prop:F estimates} asserts that $\|F(\eta_1)\|_2 = O(\varepsilon^{1/2})$ uniformly over $\eta_1 \in X_1$).
We proceed by writing
$$\NN(\eta_1+F(\eta_1)+\eta_3)+c_0^2(1-\varepsilon^2)\LL_2(\eta_1)-\KK_2(\eta_1)
=-c_0^2(1-\varepsilon^2)\NN_1(\eta_1,\eta_3)+\NN_2(\eta_1,\eta_3)+\NN_3(\eta_1,\eta_3),$$
where
\begin{align*}
\NN_1(\eta_1,\eta_3) &= \LL_2(\eta_1+F(\eta_1)+\eta_3)-\LL_2(\eta_1), \\
\NN_2(\eta_1,\eta_3) &=\KK_2(\eta_1+F(\eta_1)+\eta_3)-\KK_2(\eta_1), \\
\NN_3(\eta_1,\eta_3)&=\KK_3(\eta_1+F(\eta_1)+\eta_3)+\KK_\mathrm{r}(\eta_1+F(\eta_1)+\eta_3)\\
& \qquad \mbox{}-c_0^2(1-\varepsilon^2)\left(\LL_3(\eta_1+F(\eta_1)+\eta_3)+\LL_\mathrm{r}(\eta_1+F(\eta_1)+\eta_3)\right)
\end{align*}
and estimating these quantities using Lemma \ref{2, 3 and r estimates}.

\begin{proposition}
\label{NN1 and NN2 estimates}
The estimates
\begin{itemize}
\item[(i)]
$\|\NN_1(\eta_1,\eta_3)\|_0, \|\NN_2(\eta_1,\eta_3)\|_0\lesssim \varepsilon \nn \eta_1\nn^3+\varepsilon^{1/2} \nn \eta_1\nn^2\|\eta_3\|_2
+\varepsilon^{1/2} \nn \eta_1\nn\|\eta_3\|_2+\|\eta_3\|_2^2$,
\item[(ii)]
$\|\mathrm{d}_1\NN_1[\eta_1,\eta_3]\|_{\LL(\XX_1,L^2({\mathbb R}))},
\|\mathrm{d}_1\NN_2[\eta_1,\eta_3]\|_{\LL(\XX_1,L^2({\mathbb R}))}\lesssim \varepsilon \nn \eta_1\nn^2
+\varepsilon^{1/2} \nn \eta_1\nn\|\eta_3\|_2+\varepsilon^{1/2} \|\eta_3\|_2$,
\item[(iii)]
$\|\mathrm{d}_2\NN_1[\eta_1,\eta_3]\|_{\LL(\XX_2,L^2({\mathbb R}))},
\|\mathrm{d}_2\NN_2[\eta_1,\eta_3]\|_{\LL(\XX_2,L^2({\mathbb R}))}\lesssim \varepsilon^{1/2} \nn \eta_1\nn+\|\eta_3\|_2$
\end{itemize}
and
\begin{itemize}
\item[(iv)]
$\|\NN_3(\eta_1,\eta_3)\|_0\lesssim (\varepsilon^{1/2} \nn \eta_1\nn+\|\eta_3\|_3)^2(\nn \eta_1\nn+\|\eta_3\|_3)$,
\item[(v)]
$\|\mathrm{d}_1\NN_3[\eta_1,\eta_3]\|_{\LL(\XX_1,L^2({\mathbb R}))} \lesssim (\varepsilon^{1/2} \nn \eta_1\nn+\|\eta_3\|_3)^2$,
\item[(vi)]
$\|\mathrm{d}_2\NN_3[\eta_1,\eta_3]\|_{\LL(\XX_3,L^2({\mathbb R}))}\lesssim (\varepsilon^{1/2}\nn \eta_1\nn+\|\eta_3\|_3)(\nn \eta_1\nn+\|\eta_3\|_3)$
\end{itemize}
hold for each $\eta_1\in X_1$ and $\eta_3\in X_3$.
\end{proposition}\pagebreak

The final estimates for $\GG$ and its derivatives follow from Propositions \ref{prop:F estimates} and \ref{NN1 and NN2 estimates} by virtue of Proposition \ref{invert g}.

\begin{corollary} \label{lem:Complete GG estimates}
The estimates
\begin{itemize}
\item[(i)]
$\|\GG(\eta_1,\eta_3)\|_2\lesssim (\varepsilon^{1/2} \nn \eta_1\nn+\|\eta_3\|_2)^2(1+\nn \eta_1\nn+\|\eta_3\|_2)+\varepsilon^2\|\eta_3\|_2$,
\item[(ii)]
$\|\mathrm{d}_1\GG[\eta_1,\eta_3]\|_{\LL(\XX_1,\XX_2)} \lesssim (\varepsilon^{1/2} \nn \eta_1\nn+\|\eta_3\|_2)
(\varepsilon^{1/2}+\varepsilon^{1/2} \nn \eta_1\nn+\|\eta_3\|_2)$,
\item[(iii)]
$\|\mathrm{d}_2\GG[\eta_1,\eta_3]\|_{\LL(\XX_2,\XX_2)}\lesssim (\varepsilon^{1/2}\nn \eta_1\nn+\|\eta_3\|_2)(1+\nn \eta_1\nn+\|\eta_3\|_2)+\varepsilon^2$
\end{itemize}
hold for each $\eta_1\in X_1$ and $\eta_3\in X_3$.
\end{corollary}

\begin{theorem} \label{thm:estimate eta3}
Equation \eqref{wst eta3 eqn} has a unique solution $\eta_3 \in
X_3$ which depends smoothly upon $\eta_1 \in X_1$ and satisfies the estimates
$$
\|\eta_3(\eta_1)\|_2 \lesssim \varepsilon \nn \eta_1\nn^2, \qquad
\|\mathrm{d}\eta_3[\eta_1]\|_{\LL(\XX_1,\XX_2)} \lesssim \varepsilon \nn \eta_1\nn.
$$
\end{theorem}
{\bf Proof.}
Choosing $R_3$ and $\varepsilon$ sufficiently small and setting
$r(\eta_1)=\sigma \varepsilon \nn \eta_1 \nn^2$ for a sufficiently large value
of $\sigma>0$, one finds that
\[\|\GG(\eta_1,0)\|_2 \lesssim \tfrac{1}{2}r(\eta_1), \qquad
\|\mathrm{d}_2 \GG[\eta_1,\eta_3]\|_{\LL(\XX_2,\XX_2)}  \lesssim \varepsilon^{1/2}
\]
for $\eta_1 \in X_1$ and $\eta_3 \in \overline{B}_{r(\eta_1)}(0) \subset X_3$
(Lemma \ref{lem:Complete GG estimates}(i), (iii)).
Theorem \ref{thm:fixed-point} asserts that equation \eqref{wst eta3 eqn} has a unique solution $\eta_3$
in $\overline{B}_{r(\eta_1)}(0) \subset X_3$ which depends smoothly upon $\eta_1 \in X_1$, and
the estimate for its derivative follows from Lemma \ref{lem:Complete GG estimates}(ii).\qed\medskip

Substituting $\eta_2=F(\eta_1)+\eta_3(\eta_1)$ into \eqref{cubic eta1 eqn} yields the reduced equation
\begin{equation}
g(D)\eta_1+c_0^2\varepsilon^2K_0\eta_1+\chi(D)\left(-c_0^2(1-\varepsilon^2)\NN_1(\eta_1,\eta_3(\eta_1))+\NN_2(\eta_1,\eta_3(\eta_1))+\NN_3(\eta_1,\eta_3(\eta_1))\right)=0
\label{wst red eq v1}
\end{equation}
for $\eta_1 \in X_1$. The next step is to compute the leading-order terms in the reduced equation. To this end we write
$$\eta_1 = \eta_1^+ + \eta_1^-,$$
where $\eta_1^\pm=\chi^\pm(D)\eta_1$ and $\chi^\pm(D)$ are the characteristic functions of the sets $(\pm\omega-\delta,\pm\omega+\delta)$, so that $\eta_1^+$ satisfies the equation
\begin{equation}
g(D)\eta_1^+
+c_0^2\varepsilon^2K_0\eta_1^++\chi^+(D)\left(-c_0^2(1-\varepsilon^2)\NN_1(\eta_1,\eta_3(\eta_1))+\NN_2(\eta_1,\eta_3(\eta_1))+\NN_3(\eta_1,\eta_3(\eta_1))\right)=0
\label{Reduced equation v2}
\end{equation}
(and $\eta_1^-=\overline{\eta_1^+}$ satisfies its complex conjugate). We again begin by showing how 
Fourier-multiplier operators acting upon the function $\eta_1$ may be approximated by constants. The following result is proved in the same way as Proposition \ref{sst approximate identities}.

\begin{proposition}  \label{wst approximate identities}
The estimates
\begin{itemize}
\item[(i)]
$\partial_z \eta_1^\pm = \pm\mathrm{i}\omega\eta_1^\pm + O(\varepsilon \nn \eta_1\nn)$,
\item[(ii)]
$\partial_z^2 \eta_1^\pm =-\omega^2\eta_1^\pm+  O(\varepsilon \nn \eta_1\nn)$,
\item[(iii)]
$K_0 \eta_1^\pm = f(\omega)\eta_1^\pm + O(\varepsilon\nn \eta_1\nn)$,
\item[(iv)]
$K_0(\eta_1^+\rho_1^+) = f(2\omega)(\eta_1^+\rho_1^+) +  O(\varepsilon^{3/2}\nn \eta_1\nn \nn \rho_1 \nn)$,
 \item[(v)]
$K_0 (\eta_1^+\rho_1^-) = 2\eta_1^+ \rho_1^-+O(\varepsilon^{3/2}\nn \eta_1\nn \nn \rho_1\nn)$,
\item[(vi)]
$\FF^{-1}[g(k)^{-1}\FF[\eta_1^+\rho_1^+]]=g(2\omega)^{-1}(\eta_1^+\rho_1^+)+ O(\varepsilon^{3/2}\nn \eta_1\nn \nn \rho_1 \nn)$,
\item[(vii)]
$\FF^{-1}[g(k)^{-1}\FF[ \eta_1^+\rho_1^- ]]=g(0)^{-1}\eta_1^+\rho_1^-+O(\varepsilon^{3/2}\nn \eta_1\nn \nn \rho_1 \nn)$,
\item[(viii)]
$K_0(\eta_1^+\rho_1^+\xi_1^-) = f(\omega)(\eta_1^+\rho_1^+\xi_1^-) +O(\varepsilon^2\nn \eta_1\nn \nn \rho_1\nn \nn \xi_1 \nn)$
\end{itemize}
hold for all $\eta_1$, $\rho_1$, $\xi_1 \in \XX_1$, where the order-of-magnitude estimates are computed with respect to the $L^2({\mathbb R})$-norm.
\end{proposition}
We proceed by approximating each term in the quadratic and cubic parts of equation
\eqref{Reduced equation v2} using Corollary \ref{Formulae for 123} and Lemma \ref{wst approximate identities}.

\begin{proposition} \label{pm expansion of F}
The estimate
$$F(\eta_1) = g(2\omega)^{-1}\left(c_0^2 A(\omega)-A_0-\tfrac{1}{2}\omega^2\right)\left((\eta_1^+)^2+(\eta_1^-)^2\right) + g(0)^{-1}\left(c_0^2 B(\omega)-2A_0+\omega^2\right)\eta_1^+\eta_1^-
+\underline{O}(\varepsilon^{3/2}\nn \eta_1\nn^2),$$
where
$$A(\omega)=\tfrac{3}{2}\omega^2-\tfrac{1}{2}f(\omega)^2-f(\omega)f(2\omega)+\tfrac{1}{2}f(2\omega), \qquad
B(\omega)=\omega^2-f(\omega)^2-4f(\omega)+2,$$
holds for each $\eta_1 \in \XX_1$.
\end{proposition}
\begin{proposition} \label{Compute redeq 1}
The estimate
\begin{align*}
\chi^+(D)&\left(c_0^2(1-\varepsilon^2)\NN_1(\eta_1,\eta_3)-\NN_2(\eta_1,\eta_3)\right)\\
&=\chi^+(D)\Big(\big(2g(2\omega)^{-1}(c_0^2C(\omega)-A_0+\omega^2)(c_0^2A(\omega)-A_0-\tfrac{1}{2}\omega^2) \\
&\hspace{2cm}\mbox{}+2g(0)^{-1}(c_0^2 D(\omega)-A_0)(c_0^2B(\omega)-2A_0+\omega^2)\big)(\eta_1^+)^2\eta_1^-+ \underline{O}(\varepsilon^{3/2}\nn \eta_1 \nn^3)\Big),
\end{align*}
where
$$C(\omega)=\tfrac{3}{2}\omega^2-f(\omega)f(2\omega)+\tfrac{1}{2}f(\omega)-\tfrac{1}{2}f(\omega)^2, \qquad
D(\omega)=\tfrac{1}{2}\omega^2-\tfrac{3}{2}f(\omega)-\tfrac{1}{2}f(\omega^2),$$
holds for each $\eta_1 \in X_1$.
\end{proposition}
\begin{proposition} \label{Compute redeq 2}
The estimates
\begin{itemize}
\item[(i)]
$\chi^+(D)\KK_3(\eta_1+F(\eta_1) + \eta_3(\eta_1)) = \chi^+(D)\Big(\left(3B_0+\tfrac{1}{2}\omega^2-\tfrac{3}{2}\omega^4\right)(\eta_1^+)^2\eta_1^-+\underline{O}(\varepsilon^{3/2}\nn \eta_1 \nn^3)\Big)$,
\item[(ii)]
$\chi^+(D)\LL_3(\eta_1+F(\eta_1) + \eta_3(\eta_1)) = \chi^+(D)\Big(E(\omega)(\eta_1^+)^2\eta_1^-+\underline{O}(\varepsilon^{3/2}\nn \eta_1 \nn^3)\Big)$,
\end{itemize}
where
$$E(\omega)=2f(\omega)^2f(2\omega)-6f(\omega)\omega^2+\tfrac{13}{2}f(\omega)^2-f(\omega)f(2\omega)-4f(\omega)+\tfrac{1}{2}\omega^2,$$
hold for each $\eta_1 \in X_1$.
\end{proposition}

The higher-order terms in equation are estimated using Lemma \ref{2, 3 and r estimates}(iii).

\begin{proposition} \label{Compute redeq 3}
The estimates
\begin{itemize}
\item[(i)]
$\KK_\mathrm{r}(\eta_1+F(\eta_1) + \eta_3(\eta_1)) = \underline{O}(\varepsilon^2\nn \eta_1 \nn^4)$,
\item[(ii)]
$\LL_\mathrm{r}(\eta_1+F(\eta_1) + \eta_3(\eta_1)) = \underline{O}(\varepsilon^{3/2}\nn \eta_1 \nn^4)$
\end{itemize}
hold for each $\eta_1 \in X_1$.
\end{proposition}

\begin{corollary}
The estimate
$$
\chi^+(D)\NN_3(\eta_1,\eta_3(\eta_1)) = \chi^+(D)\Big(\left(3B_0+\tfrac{1}{2}\omega^2-\tfrac{3}{2}\omega^4-c_0^2E(\omega)\right)(\eta_1^+)^2\eta_1^- + \underline{O}(\varepsilon^{3/2}\nn \eta_1 \nn^3)\Big)
$$
holds for each $\eta_1 \in X_1$.
\end{corollary}

We conclude that the reduced equation for $\eta_1$ is the \emph{perturbed full dispersion nonlinear Schr\"{o}dinger equation}
$$g(D)\eta_1^+ + c_0^2\varepsilon^2 K_0\eta_1^+ +\chi^+(D)\Big(-4a_3
|\eta_1^+|^2\eta_1^+ +\underline{O}(\varepsilon^{3/2} \nn \eta_1 \nn^3)\Big)=0,$$
where
\begin{align*}
4a_3&= 2g(2\omega)^{-1}(c_0^2C(\omega)-A_0+\omega^2)(c_0^2A(\omega)-A_0-\tfrac{1}{2}\omega^2) \\
&\qquad\mbox{}+2g(0)^{-1}(c_0^2D(\omega)-A_0)(c_0^2B(\omega)-2A_0+\omega^2)
-3B_0-\tfrac{1}{2}\omega^2+\tfrac{3}{2}\omega^4+c_0^2E(\omega),
\end{align*}
and applying Lemma \ref{wst approximate identities}(iii), one can further simplify it to
$$g(D)\eta_1^+ + c_0^2 f(\omega) \varepsilon^2\eta_1^+ +\chi^+(D)\Big(-4a_3
|\eta_1^+|^2\eta_1^+ +\underline{O}(\varepsilon^{3/2} \nn \eta_1 \nn^3) +\underline{O}(\varepsilon^3 \nn \eta_1 \nn)\Big)=0.$$
Finally, we introduce the nonlinear Schr\"{o}dinger scaling
$$\eta_1^+(z) = \tfrac{1}{2}\varepsilon \zeta(\varepsilon z)\e^{\i \omega z},$$
noting that $I: \eta_1^+ \mapsto \zeta$ is an isomorphism $\XX_1^+:=\chi^+(D)\XX_1 \rightarrow H_\varepsilon^1({\mathbb R})$ and
$\chi^+(D)L^2({\mathbb R}) \rightarrow L_\varepsilon^2({\mathbb R})$, where\linebreak $\XX_1^+=\chi(D)\XX_1$, and choosing $R>1$ large enough so that
$\zeta_\mathrm{NLS} \in B_R(0)$ (and $\varepsilon>0$ small enough so that\linebreak $B_R(0) \subset H_\varepsilon^1({\mathbb R})$ is
contained in $I[\XX_1^+]$). We find that $\zeta \in B_R(0) \subseteq H_\varepsilon^1({\mathbb R})$ satisfies the equation
\begin{equation}
\varepsilon^{-2}g(\omega+\varepsilon D)\zeta + c_0^2 f(\omega)\zeta - a_3\chi_0(\varepsilon D)(|\zeta|^2\zeta)
+ \varepsilon^{1/2}\underline{O}^\varepsilon_0(\|\zeta\|_1)=0, \label{PFDNLS - freshly derived}
\end{equation}
which holds in $L_\varepsilon^2({\mathbb R})$.
Note that $\nn \eta_1 \nn =\varepsilon^{1/2}\|\zeta\|_1$ and the change of variable from $z$ to $Z=\varepsilon z$ introduces
an additional factor of $\varepsilon^{1/2}$ in the remainder term. Equation \eqref{wst red eq v1} is of course also invariant under the reflection $\eta_1(z) \mapsto \eta_1(-z)$,
and this invariance is inherited by \eqref{PFDNLS - freshly derived}, which is invariant under the reflection $\zeta(Z) \mapsto \overline{\zeta(-Z)}$.

\section{Solution of the reduced equation} \label{sec:existence}

In this section we find solitary-wave solutions of the reduced equations
\begin{equation}
\varepsilon^{-2}g(\varepsilon D)\zeta + 2c_0^2\zeta +2c_0^2d_0\chi_0(\varepsilon D) \zeta^2 + \varepsilon^{1/2}\underline{O}^\varepsilon_0(\| \zeta \|_1)=0,
\label{PFDKdV}
\end{equation}
and
\begin{equation}
\varepsilon^{-2}g(\omega+\varepsilon D)\zeta + c_0^2 f(\omega)\zeta - a_3\chi_0(\varepsilon D)(|\zeta|^2\zeta)
+ \varepsilon^{1/2}\underline{O}^\varepsilon_0(\|\zeta\|_1)=0. \label{PFDNLS}
\end{equation}
noting that in the formal limit $\varepsilon \rightarrow 0$ they reduce to respectively the stationary Korteweg-de Vries equation
\begin{equation}
(\tfrac{1}{8}\gamma-\tfrac{9}{8})\zeta_{ZZ}+2c_0^2\zeta+2c_0^2d_0\zeta^2=0, \label{KdV again}
\end{equation}
and the stationary nonlinear Schr\"{o}dinger equation
\begin{equation}
-a_1 \zeta_{ZZ}+a_2 \zeta -a_3|\zeta|^2\zeta =0, \label{NLS again}
\end{equation}
which have explicit (symmetric) solitary-wave solutions $\zeta_\mathrm{KdV}$ and $\pm\zeta_\mathrm{NLS}$
(equations \eqref{Explicit KdV} and \eqref{Explicit NLS}).
For this purpose we use a  perturbation argument, rewriting \eqref{PFDKdV} and
\eqref{PFDNLS} as fixed-point equations and applying the following version of the implicit-function theorem.
We again treat the cases $1 < \gamma < 9$ (`strong surface tension') and $\gamma>9$ (`weak surface tension') separately.

\begin{theorem} \label{IFT}
Let $\WW$ be a Banach space, $W_0$ and $\Lambda_0$ be open neighbourhoods of respectively $w^\star$ in $\WW$ and the origin in ${\mathbb R}$
and $\HH:  W_0 \times \Lambda_0 \rightarrow \WW$ be a function which is differentiable with respect to $w \in W_0$ for each $\lambda \in \Lambda_0$.
Furthermore, suppose that 
$\HH(w^\star,0)=0$, $\mathrm{d}_1\HH[w^\star,0]: \WW \rightarrow \WW$ is an isomorphism,
$$\lim_{w \rightarrow w^\star}\|\mathrm{d}_1\HH[w, 0]-\mathrm{d}_1\HH[w^\star,0]\|_{\LL(\WW)}=0$$
and
$$\lim_{\lambda \rightarrow 0} \|\HH(w,\lambda)-\HH(w,0)\|_{\WW}=0, \quad \lim_{\lambda \rightarrow 0} \
\|\mathrm{d}_1\HH[w,\lambda]-\mathrm{d}_1\HH[w,0]\|_{\LL(\WW)}=0$$
uniformly over $w \in X_0$.

There exist open neighbourhoods $W$ of $w^\star$ in $\WW$ and $\Lambda$ of $0$ in ${\mathbb R}$
(with $W \subseteq W_0$, $\Lambda \subseteq \Lambda_0)$ and a uniquely determined mapping
$h: \Lambda \rightarrow X$ with the properties that
\begin{itemize}
\item[(i)]
$h$ is continuous at the origin (with $h(0)=w^\star$),
\item[(ii)]
$\HH(h(\lambda),\lambda)=0$ for all $\lambda \in \Lambda$,
\item[(iii)]
$w=h(\lambda)$ whenever $(w,\lambda) \in W \times \Lambda$ satisfies $\HH(w,\lambda)=0$.
\end{itemize}
\end{theorem}

\subsection{Strong surface tension}

\begin{theorem} \label{Final strong existence thm}
For each sufficiently small value of $\varepsilon>0$ equation \eqref{PFDKdV} has a small-amplitude, symmetric
solution $\zeta_\varepsilon$ in $H_\varepsilon^1({\mathbb R})$
with $\|\zeta_\varepsilon-\zeta_\mathrm{KdV}\|_1 \rightarrow 0$ as $\varepsilon \rightarrow 0$.
\end{theorem}

The first step in the proof of Theorem \ref{Final strong existence thm} is to write \eqref{PFDKdV} as the fixed-point equation
\begin{equation}
\zeta + \varepsilon^2\big(2c_0^2\varepsilon^2 + g(\varepsilon D\big))^{-1}\left(2c_0^2d_0\chi_0(\varepsilon D)\zeta^2 + \varepsilon^{1/2}\underline{\OO}^\varepsilon_0(\| \zeta \|_1)\right)=0
\label{PFDKdV FP}
\end{equation}
for $\zeta \in H_\varepsilon^1({\mathbb R})$ and use the following elementary inequality to `replace' the nonlocal operator with a differential operator.

\begin{proposition}
The inequality
$$\left|\frac{\varepsilon^2}{2c_0^2\varepsilon^2 + g(\varepsilon k)}- \frac{1}{2c_0^2+(\frac{9}{8}-\frac{1}{9})k^2}\right| \lesssim \frac{\varepsilon}{(1+k^2)^{1/2}}$$
holds uniformly over $|k| < \delta/\varepsilon$.
\end{proposition}

Using the above proposition, one can write equation \eqref{PFDKdV FP} as
$$
\zeta+F_\varepsilon(\zeta)=0,
$$
where
$$
F_\varepsilon(\zeta)=2c_0^2d_0\left(2c_0^2-(\tfrac{9}{8}-\tfrac{1}{8}\gamma)\partial_Z^2\right)^{-1}\chi_0(\varepsilon D)\zeta^2
+ \varepsilon^{1/2}\underline{O}_1^\varepsilon(\|\zeta\|_1).
$$
It is convenient to replace this equation with
$$\zeta+\tilde{F}_\varepsilon(\zeta)=0,$$
where $\tilde{F}_\varepsilon(\zeta) = F_\varepsilon(\chi_0(\varepsilon D)\zeta)$ and study it in the fixed space $H^1({\mathbb R})$
(the solution sets of the two equations evidently coincide).
We establish Theorem \ref{Final weak existence thm} by applying Theorem \ref{IFT} with
$$\WW=H^1_\mathrm{e}({\mathbb R}) := \{u \in H^1({\mathbb R}): \mbox{$u(Z)=u(-Z)$ for all $Z \in {\mathbb R}$}\},$$
$W_0=B_R(0)$, $\Lambda_0=(-\varepsilon_0,\varepsilon_0)$ for a sufficiently small value of $\varepsilon_0$,
and
$$\HH(\zeta,\varepsilon):=\zeta+\tilde{F}_{|\varepsilon|}(\zeta)$$
(here $\varepsilon$ is replaced by $|\varepsilon|$ so that $\HH(\zeta,\varepsilon)$ is defined for $\varepsilon$ in
a full neighbourhood of the origin in ${\mathbb R}$). Observe that
$$
\HH(\zeta,\varepsilon)-\HH(\zeta,0)
=2c_0^2d_0\left(2c_0^2-(\tfrac{9}{8}-\tfrac{1}{8}\gamma)\partial_Z^2\right)^{-1}[\chi_0(|\varepsilon| D)(\chi_0(|\varepsilon| D)\zeta)^2-\zeta^2]
+ |\varepsilon|^{1/2} \underline{\OO}_1^{|\varepsilon|}(\|\zeta\|_1),
$$
and noting that
$$\lim\limits_{\varepsilon \rightarrow 0} \|\chi_0(|\varepsilon| D) - I\|_{{\mathcal L}(H^1({\mathbb R}), H^{3/4}({\mathbb R}))}=0$$
because
\begin{align*}
\| \chi_0(|\varepsilon| D) u - u \|_{3/4}^2 &= \int_{|k| > \frac{\delta}{|\varepsilon|}} (1+|k|^2)^{3/4} |\hat{u}|^2 \dk \\
& \leq \sup_{|k| > \frac{\delta}{|\varepsilon|}} (1+|k|^2)^{-1/4} \int_{|k| > \frac{\delta}{|\varepsilon|}} (1+|k|^2) |\hat{u}|^2 \dk\\
& \leq \left(1+\frac{\delta^2}{|\varepsilon|^2}\right)^{\!\!-1/4} \|u\|_1^2,
\end{align*}
that
$$\chi_0(|\varepsilon| D)(\chi_0(|\varepsilon| D)\zeta)^2-\zeta^2
=
\chi_0(|\varepsilon| D)(\chi_0(|\varepsilon| D)+I)\zeta(\chi_0(|\varepsilon| D)-I)\zeta
+(\chi_0(|\varepsilon| D)-I)\zeta^2$$
and that $H^{3/4}({\mathbb R})$ is a Banach algebra, we find that
$$\lim_{\varepsilon \rightarrow 0} \|\HH(\zeta,\varepsilon)-\HH(\zeta,0)\|_1=0, \quad \lim_{\varepsilon \rightarrow 0} \
\|\mathrm{d}_1\HH[\zeta,\varepsilon]-\mathrm{d}_1\HH[\zeta,0]\|_{\LL(H^1({\mathbb R}))}=0$$
uniformly over $\zeta \in B_R(0)$. The equation
$$\HH(\zeta,0)=\zeta+2c_0^2d_0\left(2c_0^2-(\tfrac{9}{8}-\tfrac{1}{8}\gamma)\partial_Z^2\right)^{-1}\zeta^2=0$$
has the (unique) nontrivial solution $\zeta_\mathrm{KdV}$ in $H_\mathrm{e}^1({\mathbb R})$ and it remains to show that
$$\mathrm{d}_1\HH[\zeta_\mathrm{KdV},0] =I+4c_0^2d_0\left(2c_0^2-(\tfrac{9}{8}-\tfrac{1}{8}\gamma)\partial_Z^2\right)^{-1}(\zeta_\mathrm{KdV} \cdot)$$
is an isomorphism. This result follows from the following lemma.

\begin{lemma} $ $
\begin{itemize}
\item[(i)]
The formula $\zeta \mapsto 4c_0^2d_0\left(2c_0^2-(\tfrac{9}{8}-\tfrac{1}{8}\gamma)\partial_Z^2\right)^{-1}(\zeta_\mathrm{KdV} \cdot)$ defines a compact
linear operator $H^1({\mathbb R}) \rightarrow H^1({\mathbb R})$ and $H^1_\mathrm{e}({\mathbb R}) \rightarrow H^1_\mathrm{e}({\mathbb R})$,
and in particular $\mathrm{d}_1\HH[\zeta_\mathrm{KdV},0]$ is a Fredholm operator with index $0$.
\item[(ii)]
Every bounded solution of the equation
\begin{equation}
(\tfrac{1}{8}\gamma-\tfrac{9}{8})\zeta_{ZZ}+2c_0^2\zeta+4c_0^2d_0\zeta_\mathrm{KdV}\zeta=0, \label{KdV linearised}
\end{equation}
is a multiple of $\zeta_{\mathrm{KdV},Z}$ and is therefore antisymmetric. In particular $\ker \mathrm{d}_1\HH[\zeta_\mathrm{KdV},0]$ is trivial.
\end{itemize}
\end{lemma}

Theorem \ref{KdV thm} follows from Theorem \ref{Final strong existence thm} and the following result.

\begin{proposition}
The formulae
$$\eta=\eta_1+\eta_2(\eta_1), \quad \eta_1(z)=\varepsilon^2\zeta_\varepsilon(\varepsilon z)$$
lead to the estimate
$$\eta(z) = \varepsilon^2\zeta_\mathrm{KdV}(\varepsilon z) + o(\varepsilon^2)$$
uniformly over $z \in {\mathbb R}$
\end{proposition}
{\bf Proof.} Note that
$$\|\zeta_\varepsilon - \zeta_\mathrm{KdV}\|_\infty \lesssim \|\zeta_\varepsilon - \zeta_\mathrm{KdV}\|_1 = o(1),$$
so that
$$
\eta_1(z) = \varepsilon^2 \zeta_\mathrm{KdV}(\varepsilon z)+\varepsilon^2\big(\zeta_\varepsilon(\varepsilon z) - \zeta_\mathrm{KdV}(\varepsilon z)\big)\\
= \varepsilon^2 \zeta_\mathrm{KdV}(\varepsilon z)+o(\varepsilon^2)
$$
uniformly over $z \in {\mathbb R}$. Furthermore
$$\|\eta_2(\eta_1)\|_\infty \lesssim \|\eta_2(\eta_1)\|_2 \lesssim \varepsilon^{1/2} \nn \eta_1 \nn^2 = \varepsilon^{7/2}\|\zeta_\varepsilon\|_1^2 \lesssim \varepsilon^{7/2}.\eqno{\Box}$$

\subsection{Weak surface tension}

\begin{theorem} \label{Final weak existence thm}
For each sufficiently small value of $\varepsilon>0$ equation \eqref{PFDNLS} has two small-amplitude, symmetric
solutions $\zeta^\pm_\varepsilon$ in $H^1_\varepsilon({\mathbb R})$
with $\|\zeta^\pm_\varepsilon \mp \zeta_\mathrm{NLS}\|_1 \rightarrow 0$ as $\varepsilon \rightarrow 0$.
\end{theorem}

We again begin the proof of Theorem \ref{Final weak existence thm}
by `replacing' the nonlocal operator in the fixed-point formulation
\begin{equation}\label{PFDNLS-reduced}
\zeta+\varepsilon^2\big(\varepsilon^2c_0^2 f(\omega)+g(\omega+\varepsilon D)\big)^{-1}
\left(- a_3\chi_0(\varepsilon D)(|\zeta|^2\zeta)
+ \varepsilon^{1/2}\underline{\OO}^\varepsilon_0(\|\zeta\|_1)\right)=0
\end{equation}
of equation \eqref{PFDNLS} for $\zeta \in H_\varepsilon^1({\mathbb R})$ with a differential operator.

\begin{proposition}
The inequality
$$\left|\frac{\varepsilon^2}{c_0^2 f(\omega)\varepsilon^2 + g(\omega+\varepsilon k)}- \frac{1}{a_2+a_1k^2}\right| \lesssim \frac{\varepsilon}{(1+k^2)^{1/2}}$$
holds uniformly over $|k| < \delta/\varepsilon$.
\end{proposition}

Using the above proposition, one can write equation \eqref{PFDNLS-reduced} as
$$
\zeta+\tilde{F}_\varepsilon(\zeta)=0,
$$
where
$$
\tilde{F}_\varepsilon(\zeta) = F_\varepsilon(\chi_0(\varepsilon D)\zeta), \qquad F_\varepsilon(\zeta)=-a_3\left(a_2-a_1\partial_Z^2\right)^{-1}\chi_0(\varepsilon D)(|\zeta|^2\zeta)
+ \varepsilon^{1/2} \underline{\OO}_1^\varepsilon(\|\zeta\|_1),
$$
and establish Theorem \ref{Final weak existence thm} by applying Theorem \ref{IFT} with
$$\WW=H_\mathrm{e}^1({\mathbb R}, {\mathbb C})=\{\zeta \in H^1({\mathbb R}): \mbox{$\zeta(Z) = \overline{\zeta(-Z)}$ for all $Z \in {\mathbb R}$}\},$$
$W_0=B_R(0)$, $\Lambda_0=(-\varepsilon_0,\varepsilon_0)$ for a sufficiently small value of $\varepsilon_0$
and
$$\HH(\zeta,\varepsilon):=\zeta+\tilde{F}_{|\varepsilon|}(\zeta).$$
Observe that
\begin{align*}
\HH(\zeta&,\varepsilon)-\HH(\zeta,0) \\
&=-a_3\left(a_2-a_1\partial_Z^2\right)^{-1}\Big[\chi_0(|\varepsilon| D)\big(|\chi_0(|\varepsilon| D)\zeta|^2(\chi_0(|\varepsilon| D)-I)\zeta+|\zeta|^2(\chi_0(|\varepsilon |D)-I)\zeta\\
& \hspace{2in}\mbox{}+\zeta\chi_0(|\varepsilon| D)\zeta(\chi_0(|\varepsilon| D)-I)\bar{\zeta}\big)\\
&\hspace{1.5in}\mbox{}+(\chi_0(|\varepsilon D|)-I)|\zeta|^2\zeta\Big] + |\varepsilon|^\frac{1}{2} \underline{\OO}_1^{|\varepsilon|}(\|\zeta\|_1);
\end{align*}
noting that $H^1({\mathbb R};{\mathbb C})$ is a Banach algebra, that $\chi_0(|\varepsilon| D) \rightarrow I$ in ${\mathcal L}(H^1({\mathbb R}), H^{3/4}({\mathbb R}))=0$
as $\varepsilon \rightarrow 0$ and that pointwise multiplication defines a bounded trilinear mapping $(H^1({\mathbb R}; {\mathbb C})^2 \times H^{3/4}({\mathbb R};{\mathbb C})
\rightarrow L^2({\mathbb R};{\mathbb C})$ (see H\"{o}rmander \cite[Theorem 8.3.1]{Hoermander}), one concludes that
$$\lim_{\varepsilon \rightarrow 0} \|\HH(\zeta,\varepsilon)-\HH(\zeta,0)\|_1=0, \quad \lim_{\varepsilon \rightarrow 0} \
\|\mathrm{d}_1\HH[\zeta,\varepsilon]-\mathrm{d}_1\HH[\zeta,0]\|_{\LL(H^1({\mathbb R}, {\mathbb C}))}=0$$
uniformly over $\zeta \in B_R(0)$.

The equation
$$
\HH(\zeta,0)=\zeta-a_3\left(a_2-a_1\partial_Z^2\right)^{-1}|\zeta|^2\zeta=0
$$
has (precisely two) nontrivial solutions $\pm\zeta_\mathrm{NLS}$ in $H_\mathrm{e}^1({\mathbb R}, {\mathbb C})$, which are both real, and the fact
that $\mathrm{d}_1\HH[\pm\zeta_\mathrm{NLS},0]$ is an isomorphism is conveniently established by using real
coordinates. Define $\zeta_1=\re \zeta$ and $\zeta_2=\im \zeta$, so that
$$\mathrm{d}_1\HH[\pm\zeta_\mathrm{NLS},0](\zeta_1+\mathrm{i}\zeta_2)=\HH_1(\zeta_1) + \mathrm{i} \HH_2(\zeta_2),$$
where $\HH_1: H_\mathrm{e}^1({\mathbb R}) \rightarrow H_\mathrm{e}^1({\mathbb R})$ and
$\HH_2: H_\mathrm{o}^1({\mathbb R}) \rightarrow H_\mathrm{o}^1({\mathbb R})$ are given by
$$\HH_1(\zeta_1)=\zeta_1-3a_3\left(a_2-a_1\partial_Z^2\right)^{-1}\zeta_\mathrm{NLS}^2\zeta_1, \qquad
\HH_2(\zeta_2)=\zeta_2-a_3\left(a_2-a_1\partial_Z^2\right)^{-1}\zeta_\mathrm{NLS}^2\zeta_2$$
and
\begin{align*}
H^1_\mathrm{e}({\mathbb R}) &:= \{u \in H^1({\mathbb R}): \mbox{$u(Z)=u(-Z)$ for all $Z \in {\mathbb R}$}\}, \\
H^1_\mathrm{o}({\mathbb R}) &:= \{u \in H^1({\mathbb R}): \mbox{$u(Z)=-u(-Z)$ for all $Z \in {\mathbb R}$}\}.
\end{align*}

\begin{proposition} $ $
\begin{itemize}
\item[(i)]
The formulae
$$\zeta_1 \mapsto -3a_3\left(a_2-a_1\partial_Z^2\right)^{-1}\zeta_\mathrm{NLS}^2\zeta_1, \qquad
\zeta_2 \mapsto -a_3\left(a_2-a_1\partial_Z^2\right)^{-1}\zeta_\mathrm{NLS}^2\zeta_2$$
define compact linear
operators $H^1({\mathbb R}) \rightarrow H^1({\mathbb R})$, $H_\mathrm{e}^1({\mathbb R}) \rightarrow H_\mathrm{e}^1({\mathbb R})$ and
$H_\mathrm{o}^1({\mathbb R}) \rightarrow H_\mathrm{o}^1({\mathbb R})$, and in particular $\HH_1$, $\HH_2$ are Fredholm operators
with index $0$.
\item[(ii)]
Every bounded solution of the equation
\begin{equation}
- a_1\zeta_{1ZZ} + a_2\zeta_1 - 3a_3\zeta_\mathrm{NLS}^2\zeta_1 =0 \label{NLS linearised 1}
\end{equation}
is a multiple of $\zeta_{\mathrm{NLS},Z}$ and is therefore antisymmetric, while every bounded solution of the equation
\begin{equation}
- a_1\zeta_{1ZZ} + a_2\zeta_1 - a_3\zeta_\mathrm{NLS}^2\zeta_1 =0 \label{NLS linearised 2}
\end{equation}
is a multiple of $\zeta_\mathrm{NLS}$ and is therefore symmetric. In particular $\ker \HH_1$ and $\ker \HH_2$ are trivial.
\end{itemize}
\end{proposition}

Theorem \ref{NLS thm} follows from Theorem \ref{Final weak existence thm} and the following result.

\begin{proposition}
The  formulae
$$\eta=\eta_1+F(\eta_1)+\eta_3(\eta_1), \quad \eta_1=\eta_1^++\eta_1^-, \quad \eta_1^+(z)=\tfrac{1}{2}\varepsilon\zeta_\varepsilon^\pm(\varepsilon z)\e^{\i \omega z}$$
leads to the estimate
$$\eta(z) = \pm\varepsilon\zeta_\mathrm{NLS}(\varepsilon z)\cos(\omega z) + o(\varepsilon)$$
uniformly over $z \in {\mathbb R}$.
\end{proposition}
{\bf Proof.} Note that
$$\|\zeta_\varepsilon \mp \zeta_\mathrm{NLS}\|_\infty \lesssim \|\zeta_\varepsilon \mp \zeta_\mathrm{NLS}\|_1 = o(1),$$
so that
$$
\eta_1^+(z) = \pm \tfrac{1}{2}\varepsilon \zeta_\mathrm{NLS}(\varepsilon z)\e^{\i \omega z} + \tfrac{1}{2}\varepsilon\big(\zeta_\varepsilon^\pm(\varepsilon z) \mp \zeta_\mathrm{NLS}(\varepsilon z)\big)\e^{\i \omega z} = \pm \tfrac{1}{2}\varepsilon \zeta_\mathrm{NLS}(\varepsilon z)\e^{\i \omega z}  + o(\varepsilon)
$$
uniformly over $z \in {\mathbb R}$. Furthermore
$$\|F(\eta_1)\|_\infty \lesssim \|F(\eta_1)\|_2 \lesssim \varepsilon^{1/2}\nn \eta_1 \nn^2 = \varepsilon^{3/2}\|\zeta_\varepsilon^\pm\|_1^2 \lesssim \varepsilon^{3/2}$$
and
$$\|\eta_3(\eta_1)\|_\infty \lesssim \|\eta_3(\eta_1)\|_2 \lesssim \varepsilon \nn \eta_1 \nn^2 = \varepsilon^2\|\zeta_\varepsilon^\pm\|_1^2 \lesssim \varepsilon^2.\eqno{\Box}$$

\subsection*{Acknowledgement} This work orginated with Leon Sch\"{u}tz's BSc and MSc theses at Saarland University.
The authors would like to thank Dr.\ Dan Hill for many helpful discussions concerning the radial function spaces used in Section 2.2.

\begin{appendices}
\section{Dispersion relation} \label{disprel}
In this appendix we establish the qualitative features of the dispersion relation
$$c^2 = \frac{\gamma-1+k^2}{f(k)}$$
shown in Figure \ref{Dispersion relation}.
Note that $c^2(0)=\frac{1}{2}(\gamma-1)$ and $c^2(k) \rightarrow \infty$ as $k \rightarrow \infty$. Furthermore, the calculation
\begin{equation}
\frac{\mathrm{d}c^2}{\mathrm{d}k}(k) = \frac{2kf(k)-(\gamma-1+k^2)f^\prime(k)}{f(k)^2}, \qquad f^\prime(k) = k - \frac{kI_0(k)I_2(k)}{I_1(k)^2} \label{c2 deriv}
\end{equation}
shows that
$$\frac{\mathrm{d}c^2}{\mathrm{d}k}(0)=0,$$
and it remains to determine whether $c^2$ has any critical points at positive values of $k$.

\begin{proposition}
The function
$$h(k)=1-k^2+\frac{2kf(k)}{f^\prime(k)}, \qquad k \geq 0,$$
is strictly monotone increasing.
\end{proposition}
{\bf Proof.} Observe that
\begin{equation}
h^\prime(k)=-2k+2f(k)\frac{\mathrm{d}}{\mathrm{d}k}\left(\frac{k}{f^\prime(k)}\right)+2k
=-2f(k)\big(\phi_1(k)\big)^{\!\!-2}\phi_1^\prime(k), \label{h deriv}
\end{equation}
where
$$\phi_1(k):=\frac{1}{k}f^\prime(k)= 1-\frac{I_0(k)I_2(k)}{I_1^2(k)}.$$
Barciz \cite[p.\ 257]{Baricz10} showed that for each $\nu>-1$ the function
$$\phi_\nu(k)=1-\frac{I_{\nu-1}(k)I_{\nu+1}(k)}{I_\nu^2(k)}, \qquad k\geq 0,$$
satisfies $\phi_\nu^\prime(k) < 0$ for $k>0$ with $\phi_\nu^\prime(0)=0$. It follows from equation
\eqref{h deriv} that $h^\prime(k)>0$ for $k>0$ with $h^\prime(0)=0$,
so that $h$ is strictly monotone increasing
(note that $\phi_1(k)>0$ since $\phi_1(0)=\tfrac{1}{2}$, $\phi_1(k) \rightarrow 0$ as $k \rightarrow \infty$ and $\phi_1$ is strictly monotone decreasing).\qed

Observing that $h(0)=9$ and $h(k) \rightarrow \infty$ as $k \rightarrow \infty$, we find from \eqref{c2 deriv} that for each fixed $\gamma>9$ there exists a unique $\omega>0$
with
$$\gamma=1-\omega^2 +\frac{2\omega f(\omega)}{f^\prime(\omega)}, \qquad \frac{\mathrm{d}c^2}{\mathrm{d}k}(\omega)=0,$$
while $c^2$ has no critical points at positive values of $k$ for $1<\gamma \leq 9$.
It follows that $c^2$ is a strictly monotone increasing function of $k$ for $1 < \gamma \leq 9$, while for $\gamma>9$ it has a unique local maximum at $k=0$
and a unique global minimum at $k=\omega>0$, where $\omega=h^{-1}(\gamma)>0$.

.

\section{Weakly nonlinear theory} \label{weakly nonlinear appendix}

\subsection*{Formal derivation of the KdV equation for $1<\gamma<9$}
We choose
$$c_0^2=\tfrac{1}{2}(\gamma-1),$$
write $c^2=c_0^2(1-\varepsilon^2)$ and substitute the Ansatz
$$\eta(z) = \varepsilon^2 \zeta_1(Z) + \varepsilon^4 \zeta_2(Z) + \cdots, \qquad Z=\varepsilon z,$$
into equation \eqref{GZCS}. Expanding
\begin{align*}
K_0 & = f(\varepsilon D) \\
& = \underbrace{f(0)}_{\displaystyle =2} - \tfrac{1}{2}\varepsilon^2 \underbrace{f^{\prime\prime}(0)}_{\displaystyle = \tfrac{1}{2}} \partial_Z^2 + O(\varepsilon^4),
\end{align*}
where $D=-\i\partial_Z$, we find from Corollary \ref{Formulae for 123} that
\begin{align*}
\KK_1(\eta) &= \varepsilon^2(\gamma-1) \zeta_1 + \varepsilon^4 \big(-\zeta_{1ZZ}+(\gamma-1)\zeta_2\big) + O(\varepsilon^6), \\
\KK_2(\eta) &= \varepsilon^4(-\gamma-\tfrac{1}{2}\gamma\nu^{\prime\prime}(1)+1)\zeta_1^2 + O(\varepsilon^6), \\
\\
\LL_1(\eta) &=2\varepsilon^2 \zeta_1 +\varepsilon^4(-\tfrac{1}{4}\zeta_{1ZZ}+2\zeta_2)+O(\varepsilon^6), \\
\LL_2(\eta) &=-5\varepsilon^4 \zeta_1^2 + O(\varepsilon^6)
\end{align*}
and of course $\KK_j(\eta)$, $\LL_j(\eta) = O(\varepsilon^6)$ for $j \geq 3$.

The $O(\varepsilon^2)$ component of equation \eqref{GZCS} is trivially satisfied, while the $O(\varepsilon^4)$ component 
yields the KdV equation
$$(\tfrac{1}{8}\gamma-\tfrac{9}{8})\zeta_{1ZZ}+2c_0^2\zeta_1+2c_0^2d_0\zeta_1^2=0.$$

\subsection*{Formal derivation of the NLS equation for $\gamma>9$}
We choose
$$\gamma=1-\omega^2+\frac{2\omega f(\omega)}{f^\prime(\omega)}, \qquad c_0^2 =\frac{2\omega}{f^\prime(\omega)},$$
write $c^2=c_0^2(1-\varepsilon^2)$ and substitute the Ansatz
$$\eta(z) = \varepsilon \eta_1(z,Z) + \varepsilon^2 \eta_2(z,Z) + \varepsilon^3 \eta_3(z,Z) + \cdots, \qquad Z=\varepsilon z,$$
into equation \eqref{GZCS}. Expanding
\begin{align*}
K_0 & = f(d+\varepsilon D) \\
& = f(d) - \i \varepsilon f^\prime(d)\partial_Z - \tfrac{1}{2}\varepsilon^2 f^{\prime\prime}(d) \partial_Z^2 + O(\varepsilon^3),
\end{align*}
where $d=-\i \partial_z$, $D=-\i\partial_Z$,  we find from Corollary \ref{Formulae for 123} that
\begin{align*}
\KK_1(\eta) &= \varepsilon\big((\gamma-1)\eta_1-\eta_{1zz}\big) + \varepsilon^2\big((\gamma-1)\eta_2-\eta_{2zz}-2\eta_{1zZ}\big) \\
& \qquad\mbox{}+\varepsilon^3\big((\gamma-1)\eta_3 -\eta_{3zz}-2\eta_{2zZ}-\eta_{1ZZ}\big) + O(\varepsilon^4), \\
\KK_2(\eta) & = \varepsilon^2 (A_0\eta_1^2 -\tfrac{1}{2}\eta_{1z}^2)+\varepsilon^3(2A_0\eta_1\eta_2-\eta_{1z}\eta_{1Z}-\eta_{1z}\eta_{2z})+O(\varepsilon^4), \\
\KK_3(\eta) & = \varepsilon^3 (B_0\eta_1^3 +\tfrac{1}{2}\eta_1\eta_{1z}^2+\tfrac{3}{2}\eta_{1z}^2\eta_{1zz})+O(\varepsilon^4),\\
\\
\LL_1(\eta) & = \varepsilon f(d)\eta_1 +\varepsilon^2(f(d)\eta_2-\i f^\prime(d)\eta_{1Z}) \\
& \qquad\mbox{}+\varepsilon^3(f(d)\eta_3-\i f^\prime(d)\eta_{2Z}-\tfrac{1}{2}f^{\prime\prime}(d)\eta_{1ZZ}) + O(\varepsilon^4), \\
\LL_2(\eta) & = \varepsilon^2 \big(-\tfrac{1}{2}\eta_{1z}^2-\tfrac{1}{2}(f(d)\eta_1)^2 -\eta_{1zz}\eta_1 - f(d)(\eta_1 f(d)\eta_1) + \tfrac{1}{2}f(d)\eta_1^2\big) \\
& \qquad\mbox{}+\varepsilon^3 \big(-\eta_{1z}\eta_{1Z}-\eta_{1z}\eta_{2z}-\eta_{1zz}\eta_2-\eta_1\eta_{2zz}-2\eta_1\eta_{1zZ}-(f(d)\eta_1)(f(d)\eta_2) \\
&\qquad\qquad\quad\mbox{}+\i(f(d)\eta_1)(f^\prime(d)\eta_{1Z})+\i f(d)(\eta_1 f^\prime(d)\eta_{1Z}) + \i f^\prime(d) (\eta_1 f(d)\eta_1)_Z \\
&\qquad\qquad\quad\mbox{} -f(d)(\eta_1 f(d)\eta_2) -f(d)(\eta_2 f(d)\eta_1) - \tfrac{1}{2}\i f^\prime(d)(\eta_1^2)_Z +f(d)(\eta_1\eta_2)\big)+O(\varepsilon^4), \\
\LL_3(\eta) & = \varepsilon^3 \big( -\tfrac{1}{2}(\eta_1^2\eta_{1z})_z + \tfrac{1}{2}(f(d)\eta_1)(\eta_1^2)_{zz}+(f(d)\eta_1)(f(d)(\eta_1f(d)\eta_1)) \\
& \qquad\mbox{}-\tfrac{1}{2}(f(d)\eta_1)(f(d)\eta_1^2)-(f(d)\eta_1)\eta_{1z}^2+\tfrac{1}{2}(\eta_1^2f(d)\eta_1)_{zz}+\tfrac{1}{2}f(d)(\eta_1^2\eta_{1zz}) \\
& \qquad\mbox{}-\tfrac{1}{2}f(d)(\eta_1^2f(d)\eta_1)+f(d)(\eta_1 f(d)(\eta_1 f(d)\eta_1))-\tfrac{1}{2}f(d)(\eta_1 f(d)\eta_1^2))+O(\varepsilon^4),
\end{align*}
and of course $\KK_j(\eta)$, $\LL_j(\eta) = O(\varepsilon^4)$ for $j \geq 4$.

The next step is to substitute the expressions
\begin{align*}
\eta_1(z,Z) & = \zeta_1(Z)\e^{\i\omega z} + \mbox{c.c.}, \\
\eta_2(z,Z) & = \zeta_0(Z) + \zeta_2(Z)\e^{2\i\omega z} + \mbox{c.c.}, \\
\eta_3(z,Z) & = \zeta_6(Z) + \zeta_5(Z)\e^{\i\omega z} + \zeta_4(Z)\e^{2\i\omega z} + \zeta_3(Z)\e^{3\i\omega z} + \mbox{c.c.},
\end{align*}
into the previous expansions. Noting that
$$f^{(j)}(d)(\e^{\i n \omega z}) = f^{(j)}(n \omega)\e^{\i n \omega z}, \qquad j \in {\mathbb N}_0,\ n \in {\mathbb Z},$$
we find that the $O(\varepsilon)$ component of equation \eqref{GZCS} is
$$g(\omega)\zeta_1 \e^{\i \omega z} + \mbox{c.c}=0,$$
which is satisfied because $g(\omega)=0$. The $O(\varepsilon^2)$ component yields the equation
\begin{align*}
g(0)\zeta_0 & + (2A_0-\omega^2-c_0^2 B(\omega))|\zeta_1|^2+\i g^\prime(\omega)\zeta_1^\prime\e^{\i\omega z} \\
& + \big(g(2\omega)\zeta_2 + (A_0 + \tfrac{1}{2}\omega^2-c_0^2A(\omega))\zeta_1^2\big)\e^{2\i\omega z} + \mbox{c.c.} =0,
\end{align*}
where
$$A(\omega)=\tfrac{3}{2}\omega^2-\tfrac{1}{2}f(\omega)^2-f(\omega)f(2\omega)+\tfrac{1}{2}f(2\omega), \qquad
B(\omega)=\omega^2-f(\omega)^2-4f(\omega)+2;$$
since $g^\prime(\omega)=0$ this equation is satisfied by choosing
\begin{align*}
\zeta_0 & = g(0)^{-1}(\omega^2-2A_0+c_0^2B(\omega))|\zeta_1|^2, \\
\zeta_2 & = g(2\omega)(-\omega_0 -\tfrac{1}{2}\omega^2+c_0^2 A(\omega)) \zeta_1^2.
\end{align*}

The coefficient of $\e^{\i\omega z}$ in the $O(\varepsilon^3)$ component of \eqref{GZCS} yields the equation
\begin{align}
g(\omega)\zeta_5 &+2(A_0-\omega^2-c_0^2 C(\omega))\overline{\zeta}_1\zeta_2 - \tfrac{1}{2}g^{\prime\prime}(\omega)\zeta_1^{\prime\prime}+c_0^2f(\omega)\zeta_1 \nonumber \\
&\qquad\mbox{}+ 2(A_0-c_0^2 D(\omega))\zeta_0\zeta_1 + (3B_0+\tfrac{1}{2}\omega^2-\tfrac{3}{2}\omega^4-c_0^2 E(\omega))|\zeta_1|\zeta_1|^2=0, \label{Fundamental}
\end{align}
where
\begin{align*}
C(\omega) & = \tfrac{3}{2}\omega^2-f(\omega)f(2\omega)+\tfrac{1}{2}f(\omega)-\tfrac{1}{2}f(\omega)^2, \\
D(\omega) & = \tfrac{1}{2}\omega^2 - \tfrac{3}{2}f(\omega)-\tfrac{1}{2}f(\omega)^2, \\
E(\omega) & = 2f(\omega)^2f(2\omega)-6f(\omega)\omega^2+\tfrac{13}{2}f(\omega)^2-f(\omega)f(2\omega)-4f(\omega)+\tfrac{1}{2}\omega^2.
\end{align*}
Substituting for $\zeta_0$ and $\zeta_2$ into equation \eqref{Fundamental} and setting $g(\omega)=0$ yields the nonlinear Schr\"{o}dinger equation
$$-a_1 \zeta_{1ZZ}+a_2 \zeta-a_3|\zeta_1|^2\zeta_1=0,$$
where $a_1 = \tfrac{1}{2}g^{\prime\prime}(\omega)$, $a_2 = c_0^2f(\omega)$ and
\begin{align*}
4a_3 &= 2 g(2\omega)^{-1}(c_0^2 C(\omega)-A_0+\omega^2)(c_0^2A(\omega)-A_0-\tfrac{1}{2}\omega^2) \\
& \qquad\mbox{}+2g(0)^{-1}(c_0^2D(\omega)-A_0)(c_0^2B(\omega)-2A_0+\omega^2) \\
& \qquad\mbox{}-3B_0-\tfrac{1}{2}\omega^2+\tfrac{3}{2}\omega^4+c_0^2E(\omega).
\end{align*}

\end{appendices}
\bibliographystyle{standard}
\bibliography{mdg}

\end{document}